\documentclass[12pt]{article}
\usepackage{amssymb}
\usepackage{amsmath,amsthm}
\usepackage{amsfonts}
\usepackage{graphicx}
\usepackage{graphics}
\usepackage{lineno}
\usepackage{xcolor}
\usepackage{ulem}
\usepackage[
    colorlinks=true,
    linkcolor=blue,
    citecolor=blue,
    urlcolor=blue
]{hyperref}

\numberwithin{equation}{section}
\numberwithin{footnote}{section}
\newtheorem{theorem}{Theorem}
\newtheorem{example}{Example}

\newtheorem{definition}{Definition}

\newtheorem{remark}{Remark}

\newcommand{\new}[1]{\textcolor{black}{#1}}
\newcommand{\newm}[1]{\textcolor{black}{#1}}
\newcommand{\newn}[1]{\textcolor{black}{#1}}
\newcommand{\newe}[1]{\textcolor{black}{#1}}

\begin{document}

\title{Variably Scaled RBF-RP for Scattered Data Approximation and Image Processing}

\author{
    Fatemeh Pooladi$^{a}$, 
    Maryam Mohammadi$^{a,}$\thanks{Corresponding author. E-mail: m.mohammadi@khu.ac.ir}, 
    Elisabeth Larsson$^{b}$ \\[2mm]
    \small{$^a$ \newm{Faculty} of Mathematical Sciences and Computer, Kharazmi University, Tehran, Iran.} \\
    \small{$^b$ Department of Information Technology, Uppsala University, Uppsala, Sweden.}
}

\date{\small\textsl{\today}}
\maketitle
\begin{abstract}
\new{
While radial basis function (RBF) methods are widely recognized for their approximation capabilities, their effectiveness is often limited when dealing with functions containing discontinuities, as classical RBFs are inherently designed \newm{to approximate} smooth \newm{function}s. Recent advancements, such as variably scaled radial basis functions (VS-RBFs) and variably scaled discontinuous radial basis functions (VSD-RBFs), have \newm{attempted to enhance RBF tools to interpolate functions with} finite discontinuities. However, challenges remain, particularly regarding the selection of optimal shape parameters and achieving sufficient accuracy in certain problematic scenarios. In this work, we propose an extension of these methods by incorporating radial polynomials (RPs) \newm{to address the challenges of selecting optimal shape parameters and increasing accuracy without the need to determine optimal shape parameters}, resulting in the VS-RBF-RP and VSD-RBF-RP frameworks. We present the mathematical definitions of these new bases and investigate their performance through numerical experiments. Specifically, we apply these methods to various image processing tasks to evaluate their effectiveness and observe their behavior in handling discontinuous data.
}
\vspace{.5cm}
\end{abstract}


\section{Introduction} 
RBFs have become a fundamental tool for scattered data approximation and interpolation due to their flexibility and mesh-free formulation. \newm{Many RBF methods depend on a shape parameter, whose optimal tuning remains an open problem that has received significant research attention \cite{Cavoretto21, Chen23, Esmaeilbeigi14, Fornberg08, Larsson24, Mohammadi26}.} Among the various RBFs, the multiquadric (MQ) function was first investigated by Hardy using a variable scale; however, he later preferred the constant-scale MQ because of its conceptual simplicity and computational efficiency \cite{Hardy90}. Building on Hardy's work, Kansa introduced the variably scaled MQ (VS-MQ), which depends on two parameters, and examined strategies for optimizing the scale parameter in order to improve approximation accuracy \cite{Kansa901,kansa902,kansa923,Kansa92}.

Since then, VS-RBFs have attracted considerable attention because of their potential to enhance both interpolation accuracy and numerical stability. Bozzini et al. \cite{Bozzini02} explored the relationship between VS-MQ and spline-based adaptive interpolation, while Driscoll \cite{Driscoll07} proposed algorithms for refining both the interpolation centers and their associated shape parameters. In addition, Fornberg et al. \cite{Fornberg07} studied VS-RBFs in connection with the Runge phenomenon and highlighted the difficulties involved in constructing well-conditioned variabl\newm{y} scale\newm{d} interpolation schemes.

\newm{To overcome the limitations of fixed-scale interpolation, Bozzini et al. \cite{Bozzini14} introduced a variably scaled kernel (VSK) approach that embeds d-dimensional data into a (d+1)-dimensional space, enabling the use of fixed-scale kernels to produce locally varying scalings.} Later, De Marchi et al. \cite{DeMarchi20} extended the \newm{VSK} framework to accommodate discontinuous data, leading to the development of the \newm{variably scaled discontinuous kernel (VSDK)}. More recently, several studies have investigated adaptive and data-driven strategies for selecting shape parameters, emphasizing the roles of local curvature information, learning-based prediction, and optimization techniques \cite{ Audone25, Cavoretto24, Heidari23, Sun23}. These studies highlight the continuing challenge of designing RBF\newm{s} that simultaneously achieve high accuracy, numerical stability, and computational efficiency across a wide range of applications.

Motivated by these developments, this study introduces two new bases, namely the variably scaled RBF-RP (VS-RBF-RP) and the variably scaled discontinuous RBF-RP (VSD-RBF-RP), which combine the concepts of RBF-RP, VS-RBF, and VSD-RBF.

\newm{T}his paper is organized as follows. Section 2 presents the definitions of VS-RBF and VSD-RBF. In Section 3, the RBF-RP basis is introduced, and based on the definitions provided in Section 2, the proposed VS-RBF-RP and VSD-RBF-RP bases are formulated. Section 4 \newm{presents numerical experiments and compares the proposed VS-RBF-RP and VSD-RBF-RP approaches with the classical VS-RBF and VSD-RBF methods}. Finally, Section 5 \newm{demonstrates the performance of these new bases through several image processing experiments}.

\section{Variably scaled radial basis function\newm{s}}
\new{Let $\newm{X}=\{x_1,x_2,\ldots,x_N\}\subset \Omega \subseteq \mathbb{R}^d$ denote a set of scattered data sites, and let
\begin{align*}
\newn{F}=\{\newn{f}_i=f(x_i)\}_{i=1}^N
\end{align*}
be the corresponding function values of an unknown function $f:\Omega\to\mathbb{R}$. To reconstruct $f$ by means of radial basis function (RBF) interpolation, we consider an interpolant of the form
\begin{align*}
s(x)=\sum_{i=1}^{N}\alpha_i~\phi(\varepsilon \newm{\Vert x - x_j\Vert),
  \qquad \new{x} \in \mathbb{R}^d,}
\end{align*}
where $\phi$ is a radial basis function, $\varepsilon>0$ is the associated shape parameter, and the coefficients $\alpha_i\in\mathbb{R}$ are determined by imposing the interpolation conditions
\begin{align*}
s(x_i)=\newn{f}_i,\qquad i=1,2,\cdots,N.
\end{align*}
For comprehensive treatments of RBF\newm{s} and their theoretical properties, we refer the reader to \cite{Fasshauer07, Wendland04}. \newm{T}he shape parameter, $\varepsilon$, plays a fundamental role, as it strongly influences both the approximation accuracy and the numerical stability of the resulting interpolation scheme. As noted earlier, the determination of an appropriate shape parameter based on the spatial distribution of data points, particularly through spatially varying formulations, remains a challenging and important research problem. \newm{In the following, we present the essential definitions and results required for this study}.} 
\begin{definition}\cite{Bozzini14, Romani19}\label{def1}
If a scale function $c:\new{\Omega}\rightarrow (0,\infty)$ is given, we can define a VS\newm{K} on $\new{\Omega}$ by
\begin{align}\label{eqc}
\newm{K_c( x,y)}:=\newm{K((x,c(x)),(y,c(y)))},\qquad\forall x,y\in\new{\Omega}.
\end{align}
\end{definition}
\begin{definition}\cite{DeMarchi20}(VSD\newm{K}: One dimensional case)\label{def2}
Let $\Omega =(a,b)\subset\mathbb{R}$ be an open interval, $L=\{\alpha,\beta\}$ with $\alpha,\beta\in\mathbb{R}_{>0}$, $\alpha\neq\beta$ and let $D=\{d_j,j=1,\cdots,l\}\subset\Omega$ be a set of distinct points such that $d_j <d_{j+1}$ for every $j$. Let $c:\Omega\rightarrow L$ be defined as
\newn{
\begin{align}\label{c1}
 c(x):= \begin{cases}
 \alpha, & x \in (a, d_1) \quad \text{or} \quad x \in [d_j, d_{j+1}), \quad \text{where } j \text{ is even}, \\
\beta,  & x \in [d_j, d_{j+1}), \quad \text{where } j \text{ is odd}.
\end{cases}
\end{align}
}
and
\newn{
\begin{align}\label{c12}
c(x)\vert_{[d_l,b)} := 
\begin{cases} 
\alpha, & \text{if } l \text{ is even}, \\
\beta,  & \text{if } l \text{ is odd}.
\end{cases}
\end{align}
}
With this choice of the scale function $c$ and similar to (\ref{eqc}), we call the kernel $K_c$ is a VSD\newm{K} on $\Omega$.
\end{definition}
\begin{definition}\cite{DeMarchi20}(VSD\newm{K}: Multidimensional case)\label{def3}
Let $\Omega\subset \mathbb{R}^d$ \newm{($d>1$)}, be an open subset with Lipschitz boundary, $L=\{\alpha_1,\cdots,\alpha_m\}$ real distinct values and $P = \{R_1,\cdots,R_m\}$ a partition of \newm{$\Omega$} whose elements are regions having Lipschitz boundaries. Define $c:\Omega\rightarrow L$ as
\begin{align}\label{c2}
c(x)\vert_{R_i}:=\alpha_i.
\end{align}
With this choice of the scale function and referring to (\ref{eqc}), we call the kernel $K_c$ a VSD\newm{K} on $\Omega$.
\end{definition}
\newm{Definitions \ref{def2} and \ref{def3} are special cases of the VSK formulation (Definition \ref{def1}), obtained by setting $c(x)$ to \eqref{c1}–\eqref{c12} or \eqref{c2}, respectively.}
\newm{
\begin{remark}\label{re1}
A VSK can be viewed as a fixed-parameter kernel evaluated in a higher-dimensional augmented space.
\end{remark}
\begin{remark}\label{re2}
If the kernel is an RBF, then $K_c$ in \eqref{eqc} constitutes a VS-RBF, defined as
\begin{align}\label{ker}
K_c(x,y) := \phi \left( \varepsilon r^c \right), 
\end{align}
where $r^c := \sqrt{\|x-y\|^2 + (c(x) - c(y))^2}$. Similarly, VSDK is termed a VSD-RBF if the kernel $K_c$ in Definitions \ref{def2} and \ref{def3} satisfies \eqref{ker}.
\end{remark}
}
\section{RBF-RP basis and its variably scaled extension}\label{sec2}
In this study, \new{\newm{motivated by the enhancement of RBF \newn{methods} through the addition of RPs (RBFs-RP \newn{approach}), we improve VS-RBF \newn{method} by incorporating variably scaled radial polynomials (VS-RPs), resulting in the VS-RBF-RP framework}. To this goal, we first review the definition and some properties of RBF-RP \newn{basis}. As explained in \cite{Pooladi3, Pooladi5}, when we have $N$ data points ($N>1$) there is a $m\in\mathbb{N}$ such that}
\newm{
\begin{align*}
m&=\min_{n\in\mathbb{N}}\Big(h(n)=h(0)+\sum^{n}_{t=1}(h(t)-h(t-1))\geq N\Big),\\ h(t)&= \left\{ \begin{array}{ll}
        2t+1 & \mbox{if $d=1$} ,\\
        (t+1)^2 & \mbox{if $d=2$} ,\\
        (t+1)(t+2)(2t+3)/6 & \mbox{if $d=3$}.
        \end{array} \right. 
        .
\end{align*}
Let $r_{ij} = \|x_i - x_j\|$ denote the Euclidean distance between the points $x_i$ and $x_j$, and let $r_j = (r_{1j}, \dots, r_{Nj})^T \in \mathbb{R}^N$ ($j=1, \dots, N$) be the column vector of distances from the center $x_j$. The RBF-RP basis $\langle\psi(\varepsilon r_1), \psi(\varepsilon r_2), \dots, \psi(\varepsilon r_N)\rangle$ is defined as follows}. For each center point $x_i$ ($i\in\{1,2,\cdots,N\}$), there exists \newn{a} unique RP $r_i^{2p_i}$ with \newm{$p_i\in \{0, 1,\ldots ,m\}$}, such that
\begin{align}\label{eq1}
\psi(r_i)=\phi(\varepsilon r_i)\pm\gamma r_i^{2p_i}.
\end{align} 
\newm{The parameter $\gamma$, which may depend on $\epsilon$, determines the balance between the RBF and polynomial part\newn{s}. As proposed by the authors in \cite{Pooladi3}, the RPs bases in one-, two-, and three-dimensional spaces are given as follows, respectively. The numbers above the braces indicate the number of elements corresponding to each polynomial degree:
\begin{align*}
 \{r_i^{2p_i}\}_{i\in\{1,2,\cdots ,N \}}^{p_i\in \{0,1,2,\cdots ,m\}}&= \{1, \overbrace{ r_2^2,r_3^2}^{\text{$h(1)-h(0)=2$ }}, \overbrace{ r_4^4,r_5^4}^{\text{$h(2)-h(1)=2$}},\cdots,\overbrace{ r_{N-1}^{2m},r_N^{2m}}^{\text{$N-(2m-1)$}}\},
\\
\{r_i^{2p_i}\}_{i\in\{1,2,\cdots ,N \}}^{p_i\in \{0,1,2,\cdots ,m\}}&= \{1, \overbrace{ r_2^2,\cdots ,r_4^2}^{\text{$h(1)-h(0)=3$}}, \overbrace{ r_5^4,\cdots ,r_9^4}^{\text{$h(2)-h(1)=5$}},\cdots,\overbrace{ r_{N-m^2+1}^{2m},\cdots ,r_N^{2m}}^{\text{$N-m^2$}}\},
\\
\{r_i^{2p_i}\}_{i\in\{1,2,\cdots ,N \}}^{p_i\in \{0,1,2,\cdots ,m\}}&= \{1, \overbrace{ r_2^2,\cdots ,r_5^2}^{\text{$h(1)-h(0)=4$}}, \overbrace{ r_6^2,\cdots ,r_{14}^2}^{\text{$h(2)-h(1)=9$}},\cdots,\overbrace{r_{N-\frac{m(m+1)(2m+1)}{6}+1}^{2m},\cdots ,r_N^{2m}}^{\text{$N-\frac{m(m+1)(2m+1)}{6}$}}\}.
\end{align*}
The sign $\pm$ in (\ref{eq1})} is determined by the \newn{sign in the} Taylor expansion of $\phi$ \cite{Pooladi3, Pooladi5}. For example, for Gaussian (GA) and MQ RBFs, the corresponding Taylor expansions are given by
\begin{align*}
(GA)~\phi(\new{\varepsilon} r_i)&=e^{-\varepsilon^2 r_i^2}=1-\varepsilon^2 r_i^2 +\frac{1}{2!} \varepsilon^4 r_i^4 - \cdots + (-1)^m\frac{1}{m!} \varepsilon^{2m} r_i^{2m} + \cdots,\\
(MQ)~\phi(\new{\varepsilon} r_i)&=(1+(\varepsilon r_i)^2)^{\frac{1}{2}} =1+\frac{1}{2}\varepsilon^2 r_i^2 -\frac{1}{8}\varepsilon^4 r_i^4+\cdots+ {\frac{1}{2} \choose m} \varepsilon^{2m} r_i^{2m} +\cdots.
\end{align*}
\newm{The interpolation matrix of the RBF-RP method for a given set of data points is defined as
\begin{align*}
C = A + \gamma B,
\end{align*}
where $A = [\phi(\varepsilon r_{ij})]_{i,j=1}^{N}$ denotes the standard RBF interpolation matrix, and $B = [\pm r_{ij}^{2p_{j}}]_{i,j=1}^{N}$ is the coefficient matrix corresponding to the RP component.}

An important property of the RBF-RP basis is that, due to its definition, \newn{it} behaves similarly to a standard RBF \newn{basis} when the shape parameter $\varepsilon$ is large. As explained in \cite{Pooladi1, Pooladi4}, neither \newn{the} RBF\newn{s} nor \newn{the} RP\newn{s} alone provides a suitable basis for approximating discontinuous functions, and in such cases the optimal shape parameter is typically large. Consequently, \newn{the} RBF and RBF-RP \newn{bases} behave similarly when approximating discontinuous functions.

On the other hand, VS-RBF \newn{basis} can be interpreted as a special case of an RBF defined in a higher-dimensional space. Therefore, following the same strategy used in the RBF-RP construction, we incorporate \new{VS-}RP into \newn{the} VS-RBF and VSD-RBF \newn{bases} in order to improve the accuracy and stability of these bases. The resulting bases are referred to as VS-RBF-RP and VSD-RBF-RP. 

To define VS-RBF-RP and VSD-RBF-RP \newn{bases}, we first use the definitions of VS-RBF and VSD-RBF given in Definitions \ref{def1}–\ref{def3}. Then we explain how the parameter $\gamma$ is defined when the shape parameter is variable rather than constant. Finally, we define the RP corresponding to a variable \newm{scale (VS-RP)}.

\newm{Recently, the authors in \cite{Pooladi5} established a geometric framework to define the weighting parameter $\gamma$ in the RBF-RP method. In this approach, the flatness of an RBF is evaluated by examining the proximity of its curvature at the origin to zero. Based on this characteristic, the following functional relation was introduced:
\begin{align*}
\gamma = u\left(\left.\kappa_\phi(r)\right|_{r=0}\right),
\end{align*}
where $\left.\kappa_\phi(r)\right|_{r=0}$ denotes the curvature of $\phi(\varepsilon r)$ at the origin, and $u$ is a continuous, strictly decreasing function satisfying:
\begin{align*}
0 < u(t) \le 1, \quad \lim_{t \to 0} u(t) = 1, \quad \text{and} \quad \lim_{t \to \infty} u(t) = 0.
\end{align*}
This monotonicity guarantees that as the RBF \newn{profile} becomes more peaked (i.e., having a higher curvature), the contribution of the polynomial term systematically diminishes.
\begin{remark}\cite{Heidari23, Pooladi5}\label{re4}
For standard RBF \newn{bases}, including GA, MQ, Inverse Multiquadric (IMQ), and Inverse Quadratic (IQ), the curvature at the origin in any spatial dimension is given by:
\begin{align*}
\left.\kappa_\phi(r)\right|_{r=0} = \eta \varepsilon^2,
\end{align*}
where $\eta$ equals $2$, $1$, $1$, and $2$, respectively.
\end{remark}
Furthermore, they proposed the decay function $u(t) = \exp(-\sigma t)$ with $\sigma > 0$, leading to 
\begin{align}\label{eq2}
\gamma(\varepsilon) = \exp(-\sigma \eta \varepsilon^2).
\end{align}
In particular, selecting $\sigma = 1/\eta$ neutralizes the dependency on the specific RBF type, thereby yielding an RBF-RP formulation that is independent of the choice of basis function.}
\newm{
To generalize this curvature-based formulation of $\gamma$ to variably scaled methods (VS-RBF-RP and VSD-RBF-RP), we first analyze the curvature behavior of VS-RBF and VSD-RBF \newn{bases} in the following.
}
\begin{theorem}\label{theo1}
The curvature of the \newe{VS-RBF} \newn{at the origin} is equal to \newn{the curvature of the} \newe{RBF}.
\end{theorem}
\begin{proof}
According to Remark \ref{re1}, the result follows directly.
\end{proof}
\newm{Clearly, the curvature of the VSD-RBF \new{basis function} is \newn{identical} to that of VS-RBF \newn{basis function}.} According to (\ref{eq2})\newm{, Remark \ref{re4}} and Theorem \ref{theo1}, the parameter $\gamma$ for VS-RBF-RP and VSD-RBF-RP \newn{methods} is therefore given by
\begin{align*}
\gamma = \new{\mathrm{e}^{-\varepsilon^2}}.
\end{align*}
A\newm{ccording to Remark \ref{re2}, if kernel} is an RBF, the interpolation matrix elements associated with the VS-RBF basis are given by
\begin{align*}
\phi(\newm{\varepsilon \sqrt{\|x_j-x_i\|^2+(c(x_j)-c(x_i))^2}}).
\end{align*}
Based on the formulation of the RBF-RP method in \eqref{eq1}, we now introduce the general variabl\newm{y} scaled RBF with radial polynomials (VS-RBF-RP) as follows.
\begin{definition}(VS-RBF-RP)\label{def4}
\new{Let $c: \newm{\Omega} \rightarrow (0,\infty)$ be a given scale function and let $\phi$ denote an infinitely smooth radial basis function. For any points $x_i, x_j \in \newm{\Omega}$, the VS-RBF-RP \newn{basis} on $\newm{\Omega}$ is defined as}
\begin{align}\label{eq5}
K_c(x_{\newm{i}},x_{\newm{j}}) :&= \phi(\newm{\varepsilon\sqrt{\|x_i-x_j\|^2+(c(x_i)-c(x_j))^2}})\pm \new{\mathrm{e}^{-\varepsilon^2}}(\|x_{\newm{i}}-x_{\newm{j}}\|^2+(c(x_{\newm{i}})-c(x_{\newm{j}}))^2)^{p_{\newm{j}}},\nonumber\\
&= \phi(\newm{\varepsilon \sqrt{r_{ij}^2+ c_{ij}^2}}) \pm \mathrm{e}^{-\varepsilon^2}(r_{ij}^2+c_{ij}^2)^{p_{\newm{j}}},
\end{align}
\new{where $c_{ij} = c(x_{\newm{i}})-c(x_{\newm{j}})$, and the exponent $p_{\newm{j}}$ represents the degree of the radial polynomial for the $(d+1)$-dimensional basis.}
\end{definition}
\new{As a specialized variant, when the scale function exhibits discontinuities across subdomains to capture sharp gradients or internal interfaces, the variabl\newm{y} scaled discontinuous RBF with radial polynomials (VSD-RBF-RP \newn{method}) is established.
\begin{definition}(VSD-RBF-RP)\label{def5}
\newm{Following Remark \ref{re2}, b}y adopting the scale function $c(x)$ defined in Definition~\ref{def2} (for the one-dimensional case) or Definition~\ref{def3} (for the multidimensional case), and incorporating it into \eqref{eq5}, the resulting kernel $K_c$ is called a VSD-RBF-RP \newn{basis} on $\Omega$.
\end{definition}
}
\newm{The interpolation matrix of the VS-RBF-RP method for a given set of data points is defined as
\begin{align*}
C^c = A^c + \gamma B^c,
\end{align*}
where $A^c = [\phi(\varepsilon r^c_{ij})]_{i,j=1}^{N}$ denotes the standard VS-RBF interpolation matrix, and $B^c = [(r^c_{ij})^{2p_j}]_{i,j=1}^{N}$ represents the coefficient matrix corresponding to the VS-RP component. Here,
\begin{align*}
r^c_{ij} = \sqrt{\|x_i - x_j\|^2 + (c(x_i) - c(x_j))^2}
\end{align*}
defines the augmented Euclidean distance between points $x_i$ and $x_j$. Note that in the VSD-RBF-RP framework, the coefficient matrix is constructed based on the spatial scaling function $c(x)$.}
\begin{remark}\label{re5}
In the numerical experiments reported by \cite{DeMarchi20}, the shape parameter was fixed at $\varepsilon = 1$. Our numerical findings, presented in Section~\ref{sec4}, confirm that \newm{the accuracy achieved with} this choice is \newm{very close to that obtained with the optimal shape parameter} for the test problems under consideration. \newm{This sensitivity analysis indicates that for the standard VSD-RBF method, choosing a shape parameter other than $\varepsilon = 1$ may lead to a degradation in performance.} Moreover, through a comprehensive parametric study across a wider range of $\varepsilon$, we demonstrate that the proposed VSD-RBF-RP method consistently outperforms the standard VSD-RBF approach across varying shape parameter values.
\end{remark}
\section{Numerical experiments and comparison with classical VS-RBF and VSD-RBF methods}\label{sec4}
In this section, the proposed VS-RBF-RP and VSD-RBF-RP methods are compared with previously developed approaches. To this end, several test problems from \cite{Bozzini14, DeMarchi20} are considered. \newn{\newn{Solely for the purpose of benchmark comparisons and demonstrating the best-case performance of each method, a} near-optimal shape parameter $\varepsilon$ is determined by minimizing the RMSE over $\varepsilon \in [0.01, 100]$ \newn{on a logarithmically spaced grid with a step size of $10^{0.1}$. We emphasize that such an exhaustive search is employed strictly for validation, as it becomes computationally prohibitive in practical applications where the exact solution is unknown.}}
\subsection{A continuous case}
The original motivation for introducing the VS-RBF method was to enhance numerical stability. In \cite{Bozzini14}, the authors investigated the condition numbers of interpolation matrices associated with different point distributions for the Runge function \newm{$f(x)=\dfrac{1}{1 + 25x^2}$}. Tables \ref{table:1} and \ref{table:2} compare the results of \newm{interpolating this function using different bases: RBF, VS-RBF, and VS-RBF-RP} for Chebyshev and Halton points on the interval $[-1,1]$.

\begin{table}[h!]
\centering
\begin{tabular}{lccc}
\hline
\newn{Basis} & Condition \newm{number} & \newm{RMSE}\\
\hline
\newm{RBF}& $3.549\times 10^{17}$ & $3.203\times 10^{-6}$  \\
\newm{RBF-RP}& $3.549\times 10^{17}$ & $3.203\times 10^{-6}$  \\
\newm{VS-RBF}& $4.563\times 10^{4}$ & $7.527\times 10^{-6}$ \\
\newm{VS-RBF-RP}& $4.563\times 10^{4}$ & $7.527\times 10^{-6}$  \\
\hline
\end{tabular}
\caption{\newm{Comparison of RBF, RBF-RP, VS-RBF, and VS-RBF-RP \newn{bases} (using G\newn{A} basis) with Chebyshev points ($N=55$).}}
\label{table:1}
\end{table}

\begin{table}[h!]
\centering
\begin{tabular}{lccc}
\hline
\newn{Basis} & Condition \newm{number} & \newm{RMSE}\\
\hline
\newm{RBF}& $1.015\times 10^{13}$ & $5.810\times 10^{-6}$ \\
\newm{RBF-RP}& $1.015\times 10^{13}$ & $5.810\times 10^{-6}$  \\
\newm{VS-RBF}& \newn{$6.108$} & \newn{$3.903\times 10^{-3}$}  \\
\newm{VS-RBF-RP}& \newn{$6.108$} & \newn{$3.903\times 10^{-3}$} \\
\hline
\end{tabular}
\caption{\newm{Comparison of RBF, RBF-RP, VS-RBF, and VS-RBF-RP \newn{bases} (using G\newn{A} basis) with Halton points ($N=55$).}}
\label{table:2}
\end{table}
 As described in \cite{Bozzini14}, the scale function $c(x)$ is defined as 
\begin{align*}
\sqrt{1 - x^{2}}
\end{align*}
for Chebyshev points\newn{. Furthermore, for the Halton points, it is defined by the recursive rule}
\begin{align*}
c_{j+1}=c_j+\sqrt{4h\newm{_D}^2-(x_{j+1}-x_j)^2},
\end{align*}
\newn{Let $X = \{x_1, x_2, \dots, x_n\}$ be the sorted set of 1D Halton points in the interval $[-1, 1]$, such that $-1 \le x_1 < x_2 < \dots < x_n \le 1$. The fill distance $h$ is defined as:
\begin{align*}
h_D = \max \left\{ x_1 - (-1), \max_{1 \le i < n} \left( \frac{x_{i+1} - x_i}{2} \right), 1 - x_n \right\}.
\end{align*}
}
All experiments were conducted using 55 data points, with the \newm{fixed} shape parameter $\varepsilon = 7.943$. As observed in Table \ref{table:1}, the results obtained from VS-RBF and VS-RBF-RP \newn{bases} are identical \newn{for both point distributions.} 

\newn{A clear trade-off between numerical stability and approximation accuracy is observed between the two lifting strategies. From a stability perspective, the recursive Halton lifting significantly mitigates matrix ill-conditioning, achieving an exceptional matrix condition number ($\approx 6.108$) compared to the Chebyshev lifting ($\approx 4.563 \times 10^4$). Conversely, the Chebyshev nodes yield a higher approximation accuracy, achieving an $\text{RMSE}$ of $7.527 \times 10^{-6}$ for the VS-RBF formulations, whereas the Halton nodes result in an $\text{RMSE}$ of $\approx 3.903 \times 10^{-3}$. This behavior highlights a classical trade-off in meshless methods: while the Halton-based discrete manifold maximizes solver robustness and alleviates matrix degeneracy, the smooth Chebyshev lifting preserves superior boundary resolution and interpolation accuracy.}

\subsection{Discontinuous cases}
In this subsection, \new{we evaluate} the accuracy and convergence \new{properties} of the VSD-RBF-RP basis \new{in comparison} with the VSD-RBF method. \new{To this end, t}wo benchmark tests from \cite{DeMarchi20} are \new{conducted}. In both \new{test cases, $N=625$} clustered Halton points \new{distributed} in $[-1,1]^2$ \new{serve as the centers, with a total of $M=2500$} evaluation points. \new{To ensure a fair comparison, the same scaling functions, $c_1$ and $c_2$, as employed in \cite{DeMarchi20} are adopted. Furthermore, following the parameter settings in \cite{DeMarchi20}, \newm{the results for the condition number are reported for this specific value ($\varepsilon = 1$)}.} 
\begin{example}\label{exa1}
Let the test function $f_1$ and the associated scale function $c_1$ be defined by
\begin{align*}
f_1(x,y)&=\begin{cases}
 \mathrm{e}^{-(x^2+y^2)}, & x^2+y^2\leq 0.6,\\
x+y, & x^2+y^2>0.6,
\end{cases}\\
c_1(x,y)&=\begin{cases}
1, & x^2+y^2\leq 0.6,\\
2, & x^2+y^2>0.6,
\end{cases}
\end{align*}
Figure \ref{fig1} (left) shows that the VSD-RBF-RP basis achieves a visibly higher accuracy compared with the other bases. Moreover, Figure \ref{fig1} (right) illustrates the convergence behavior of the new basis.
\begin{figure}
	\centering
        \setkeys{Gin}{width=7cm}
        \includegraphics{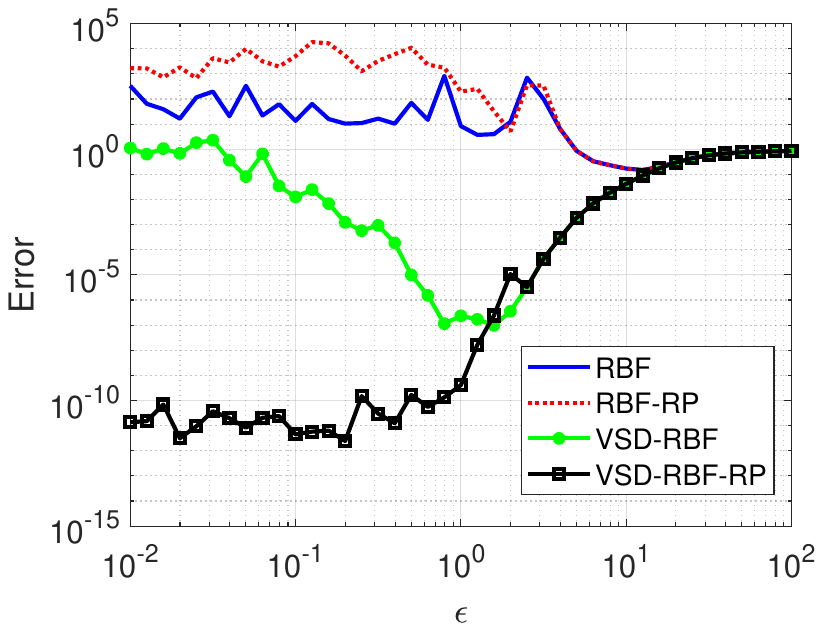}
        \includegraphics{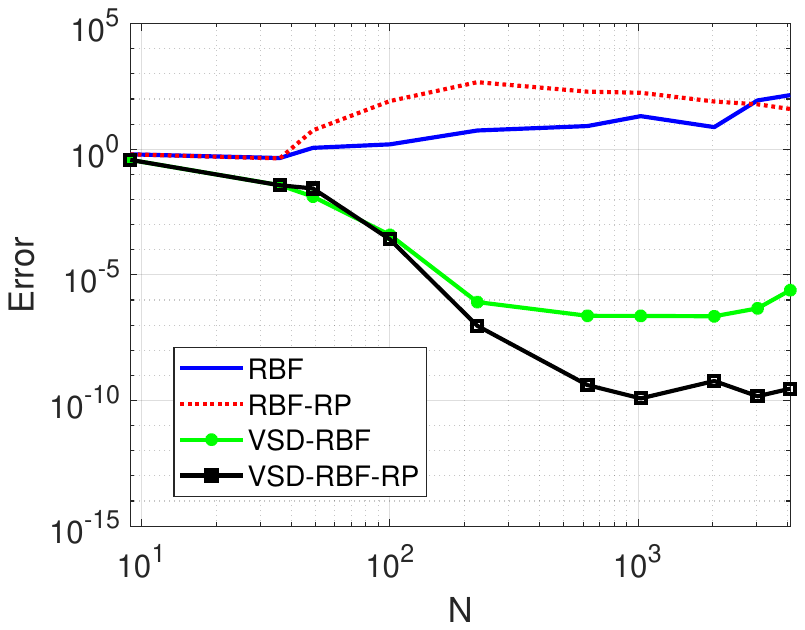}
	  \caption{\newm{RMSE comparison for $f_1$ across various shape parameter values with $N = 625$ and $M = 2500$; (left); convergence behavior with respect to the number of centers $N$ at a fixed $\varepsilon = 1$ (right).}}
\label{fig1}
\end{figure}

\newm{Figures~\ref{fig2} and~\ref{fig3} display the exact test function $f_1$ and its corresponding approximate solutions obtained using different basis functions and shape parameters, respectively. As can be observed, the VSD-RBF-RP method proves to be a superior choice across the entire range of the shape parameter.}
\begin{figure}
	\centering
        \setkeys{Gin}{width=8cm}
        \includegraphics{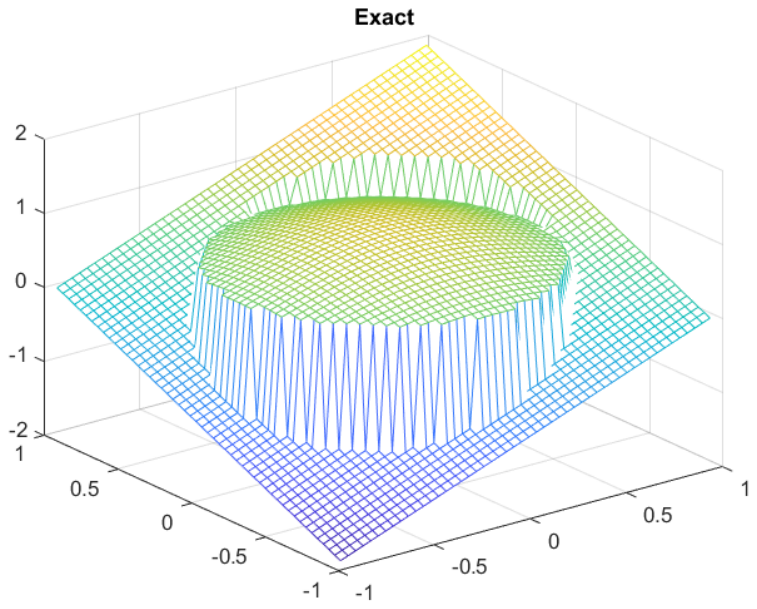}
	  \caption{\newm{Exact test function $f_1$.}}
\label{fig2}
\end{figure}
\begin{figure}
	\centering
        \setkeys{Gin}{width=4.4cm,height=4.5cm,keepaspectratio}
        \includegraphics{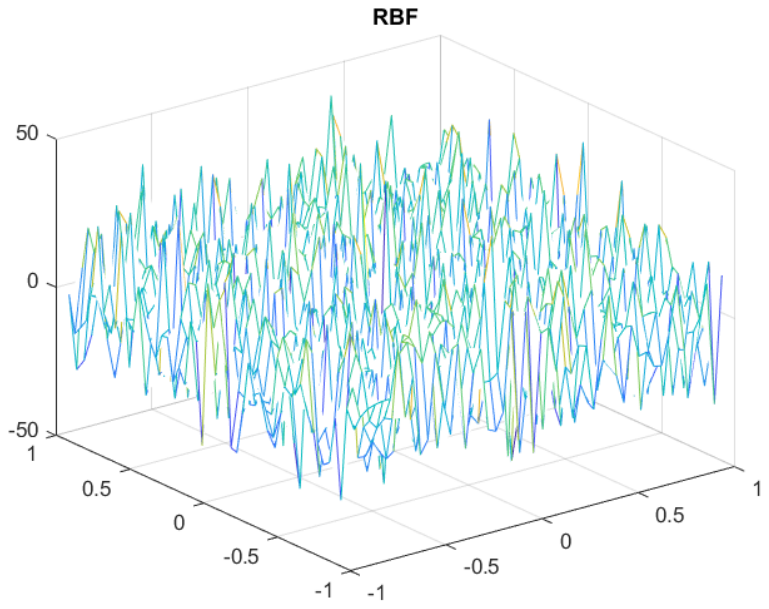}
        \includegraphics{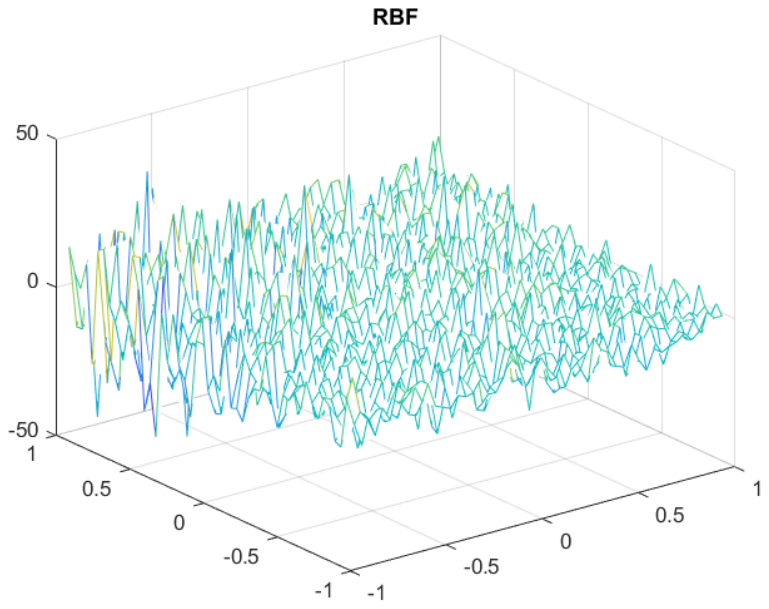}
        \includegraphics{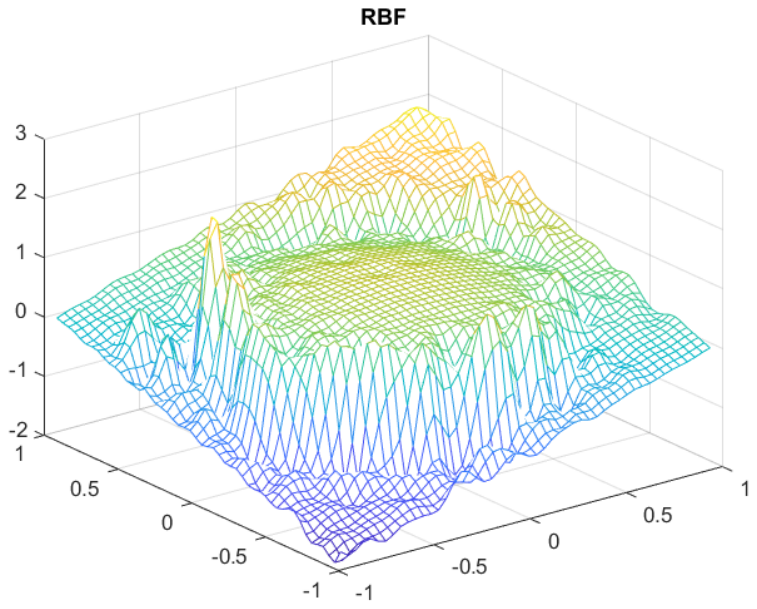}
        
        \includegraphics{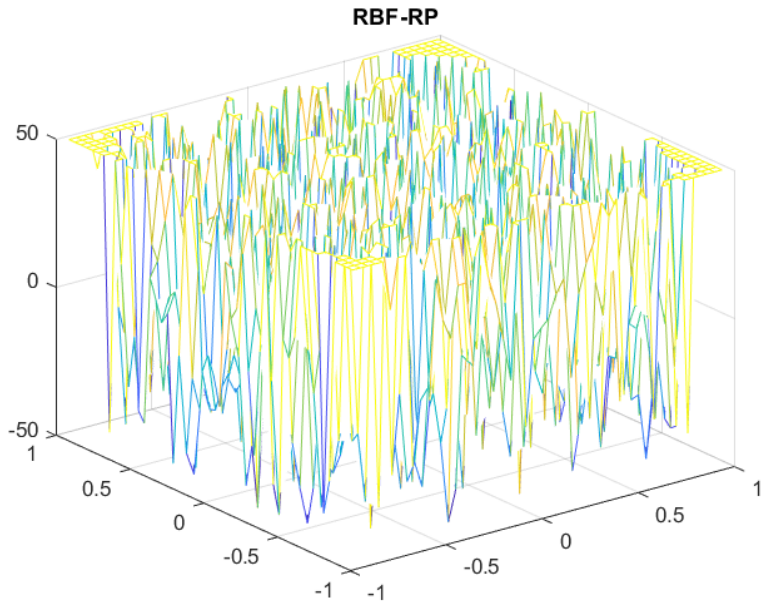}
        \includegraphics{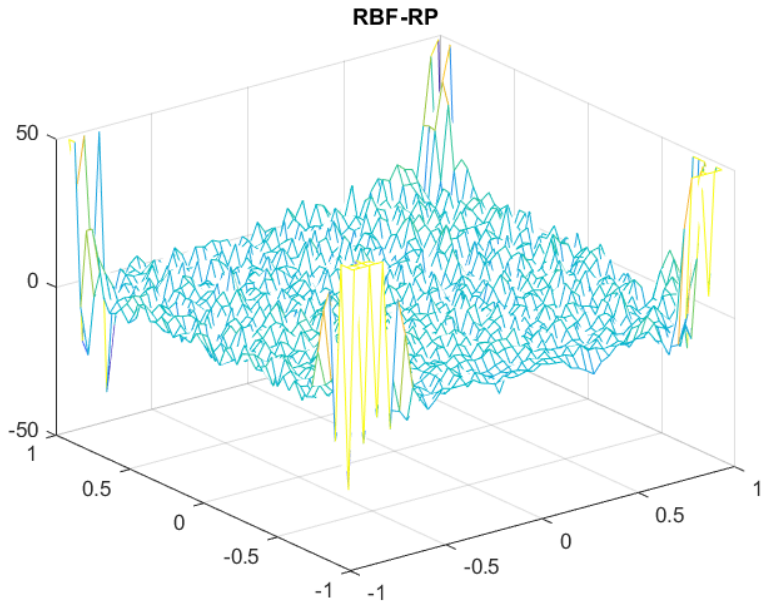}
        \includegraphics{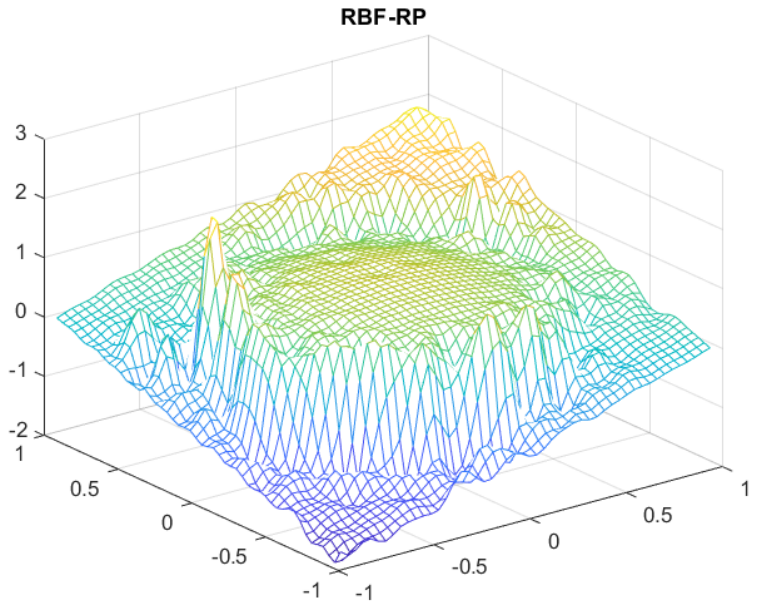}
        
        \includegraphics{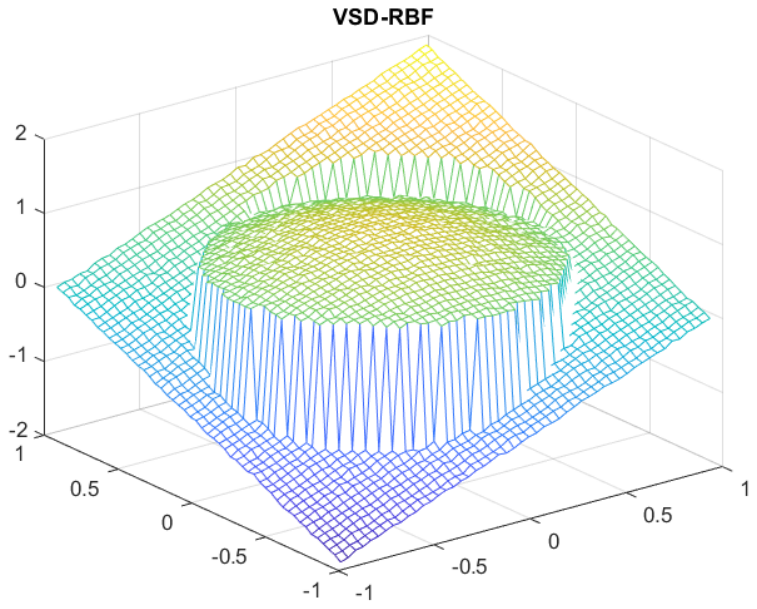}
        \includegraphics{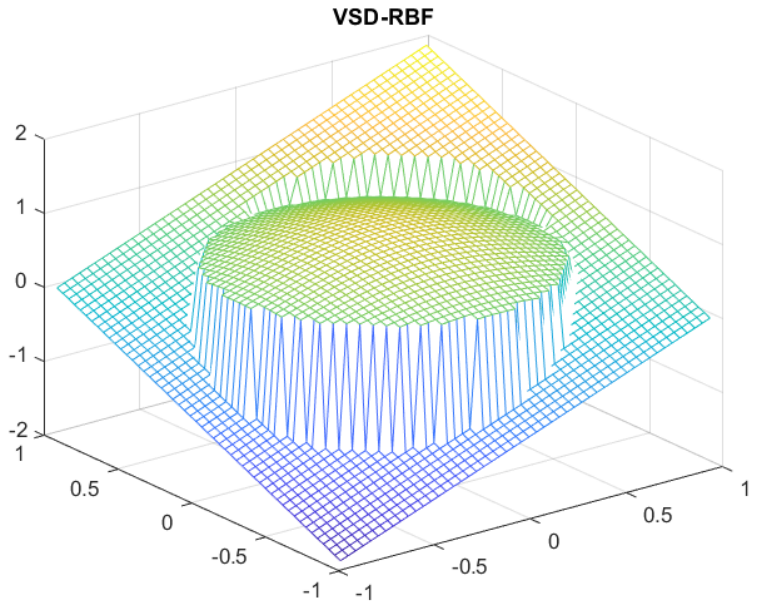}
        \includegraphics{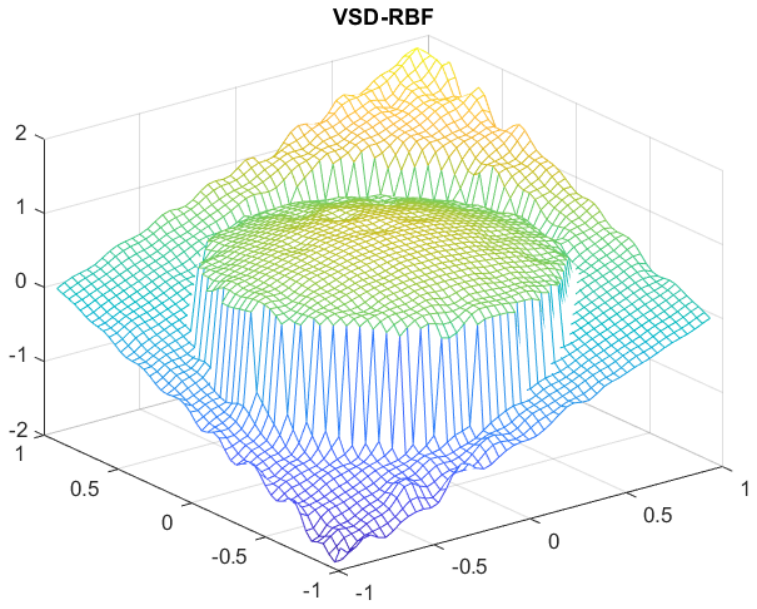}
        
        \includegraphics{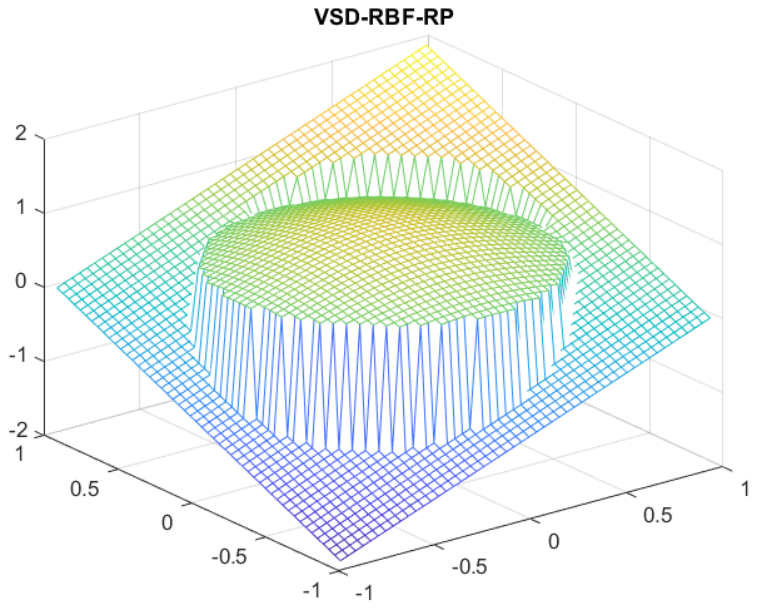}
        \includegraphics{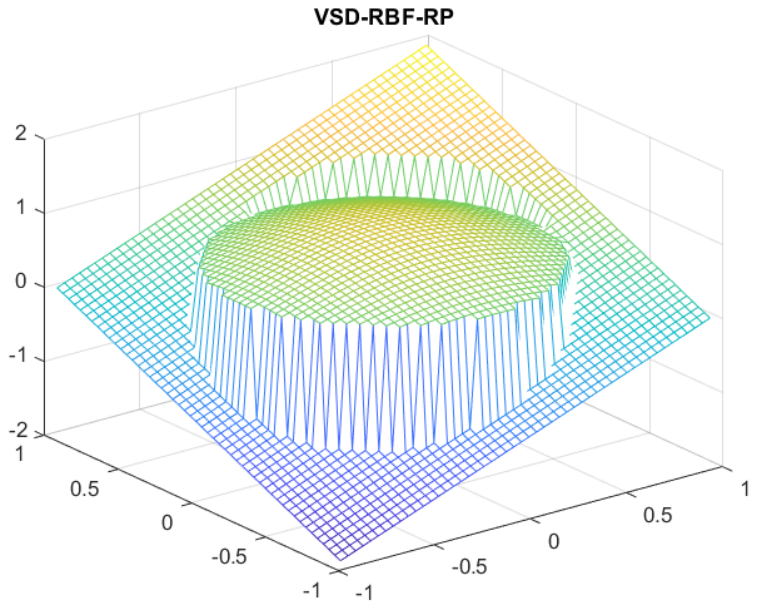}
        \includegraphics{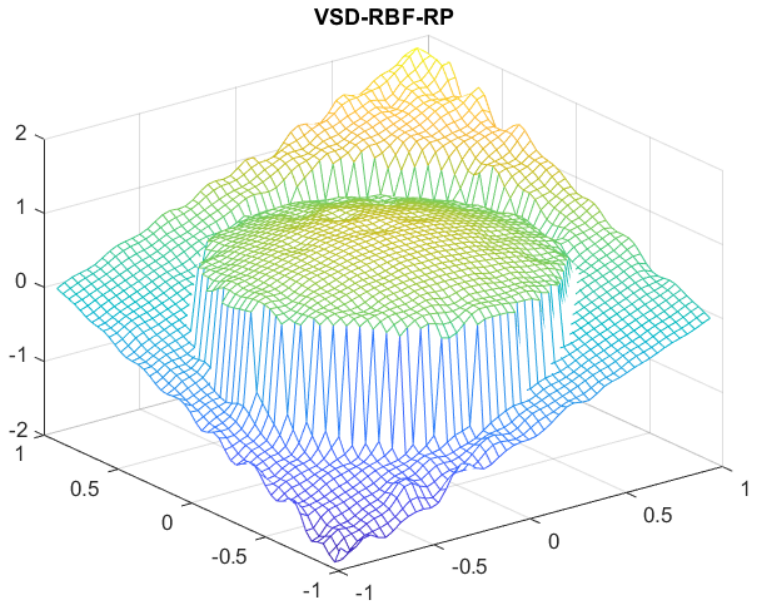}
	  \caption{\newm{Comparison of the approximate solutions for $f_1$ with $N = 625$ and $M = 2500$. Rows (top to bottom): RBF, RBF-RP, VSD-RBF, and VSD-RBF-RP \newn{methods}. Columns (left to right): $\varepsilon=0.1,1,$ and $10$.}}
\label{fig3}
\end{figure}
\end{example}
\begin{example}\label{exa2}
Next, we examine the accuracy and convergence of the bases for the function $f_2$, defined together with its scale function $c_2$ by
\begin{align*}
f_2(x,y)&=\begin{cases}
2(1-\mathrm{e}^{-(y+0.5)^2}), & \vert x\vert\leq 0.5, \vert y\vert\leq 0.5, \\
4(x+0.8), & -0.8\leq x\leq -0.65, \vert y\vert\leq 0.8, \\
0.5, & 0.65\leq x\leq 0.8, \vert y\vert\leq 0.2,  \\
0, & \text{otherwise}
\end{cases}\\
c_2(x,y)&=\begin{cases}
1, & \vert x\vert\leq 0.5, \vert y\vert\leq 0.5, \\
2, & -0.8\leq x\leq -0.65, \vert y\vert\leq 0.8, \\
3, & 0.65\leq x\leq 0.8, \vert y\vert\leq 0.2,  \\
4, & \text{otherwise}
\end{cases}
\end{align*}

Compared to $f_1(x,y)$, the test function $f_2(x,y)$ presents a significantly more challenging approximation problem \newm{due to a localized rectangular subdomain, specifically $0.65 \leq x \leq 0.8$ and $|y| \leq 0.2$, where the function abruptly transitions to a constant value of $0.5$ (see Figure~\ref{fig4} for an illustration of $f_2$). This sharp feature creates discontinuities at its boundaries that cannot be effectively captured by global approximation bases, including the polynomial term in the RBF-RP framework. Consequently, the attempt to approximate this localized abrupt change induces severe Runge- and Gibbs-type oscillations, which degrade the global approximation accuracy and explain why the VSD-RP fails to yield any noticeable improvement for $f_2(x,y)$.
}
\begin{figure}
	\centering
        \setkeys{Gin}{width=8cm}
        \includegraphics{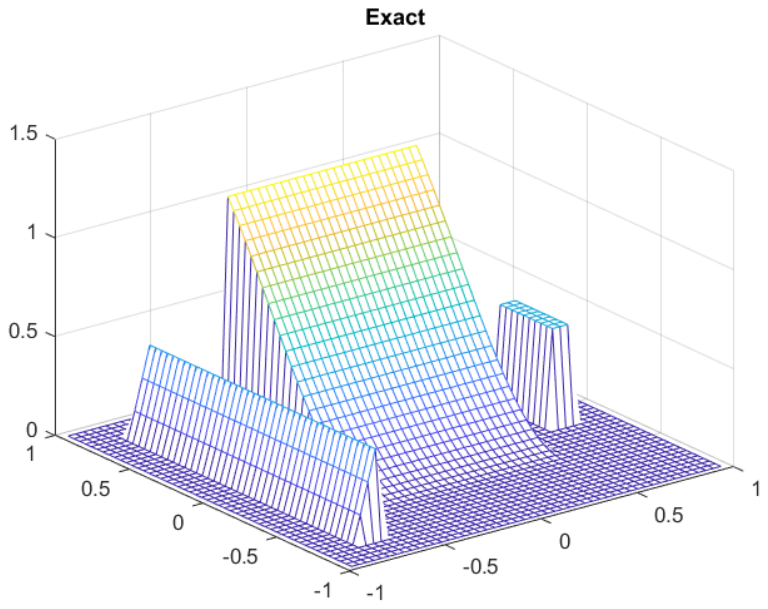}
	  \caption{\newm{Exact test function $f_2$ with $N = 625$ and $M = 2500$.}}
\label{fig4}
\end{figure}

\new{This theoretical challenge is directly corroborated by our numerical results.} As shown in Figure~\ref{fig5}, the VSD-RBF-RP method initially fails to show an advantage over VSD-RBF \newn{method}. \newm{To investigate the underlying behavior of VSD-RBF-RP \newn{method} in this setting, Figure~\ref{fig6} illustrates the reconstructed surfaces of $f_2(x,y)$ using the standard RBF, RBF-RP, VSD-RBF, and VSD-RBF-RP methods. As anticipated, regardless of the shape parameter value, the error generated by VSD-RBF-RP \newn{approach} is markedly higher within the localized region $0.65 \leq x \leq 0.8$, $|y| \leq 0.2$ compared to the rest of the computational domain.}

\begin{figure}
	\centering
        \setkeys{Gin}{width=7cm}
        \includegraphics{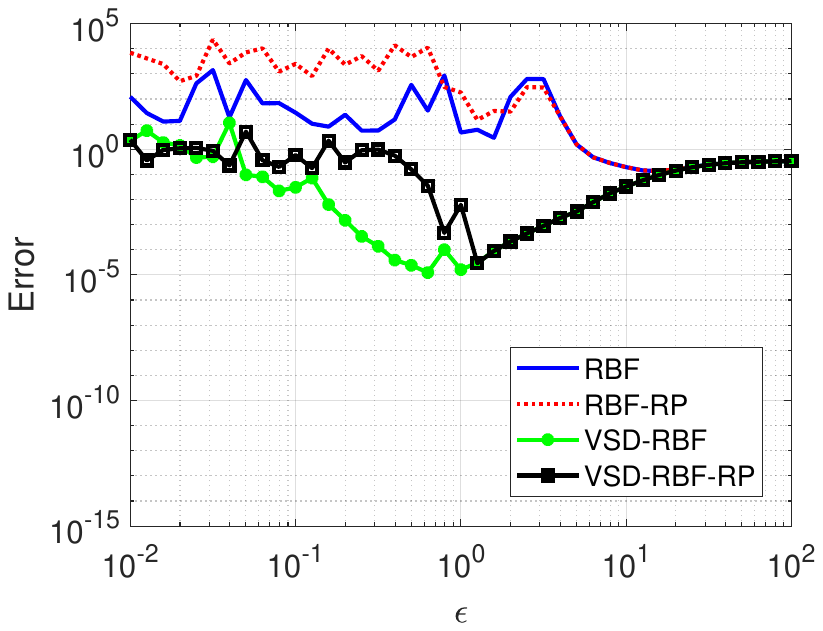}
	  \caption{\newm{RMSE comparison for $f_2$ across various shape parameter values with $N = 625$ and $M = 2500$.}}
\label{fig5}
\end{figure}
\begin{figure}
	\centering
        \setkeys{Gin}{width=4.4cm,height=4.5cm,keepaspectratio}
        \includegraphics{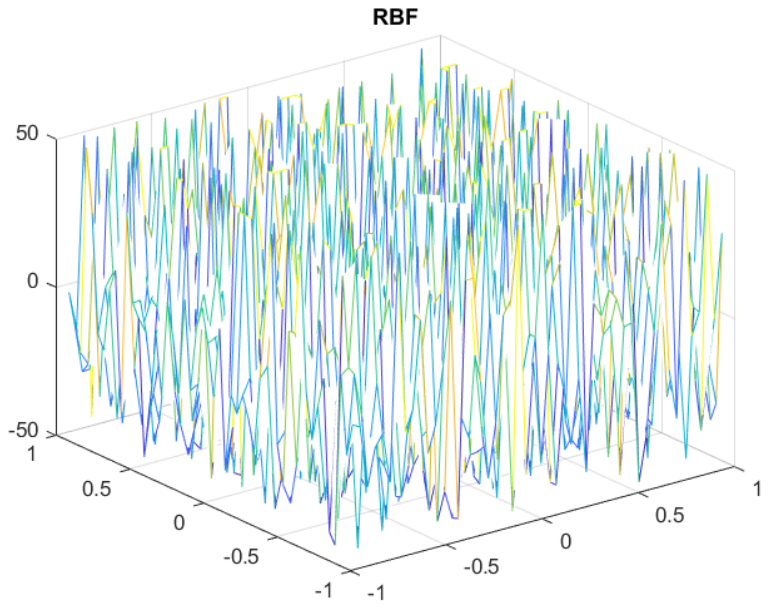}
        \includegraphics{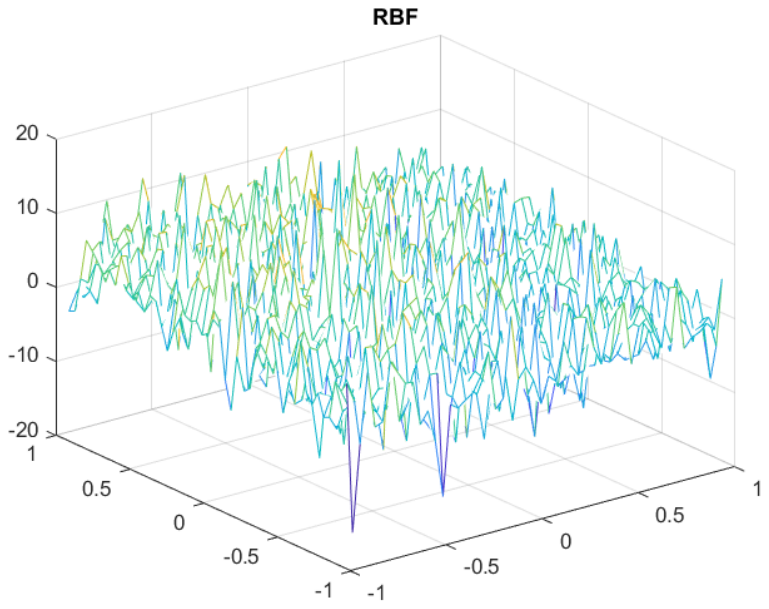}
        \includegraphics{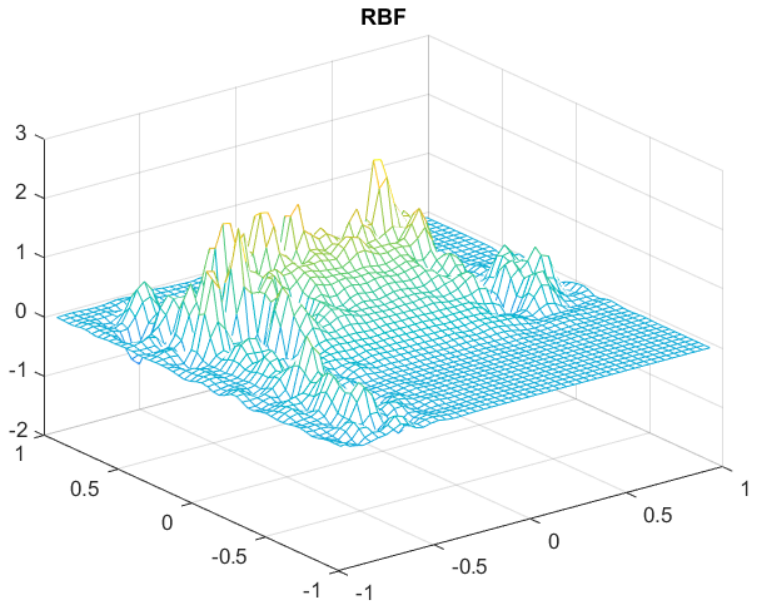}
        
        \includegraphics{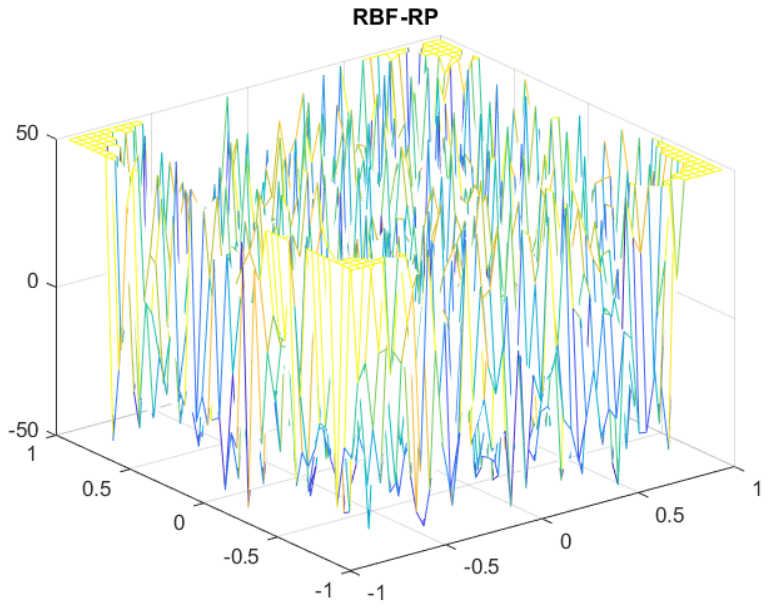}
        \includegraphics{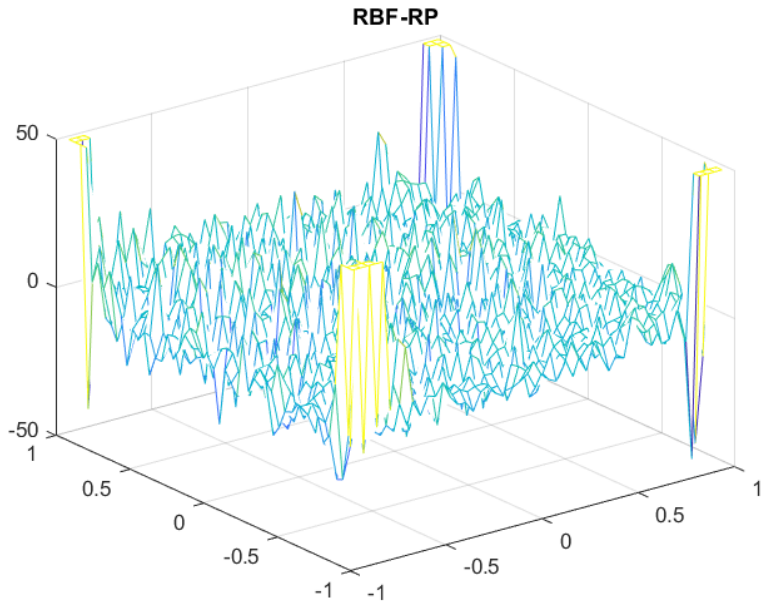}
        \includegraphics{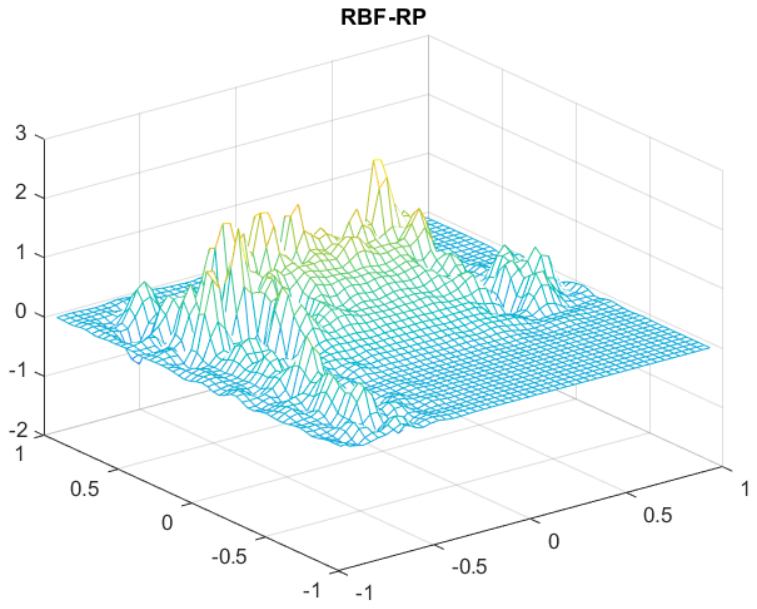}
        
        \includegraphics{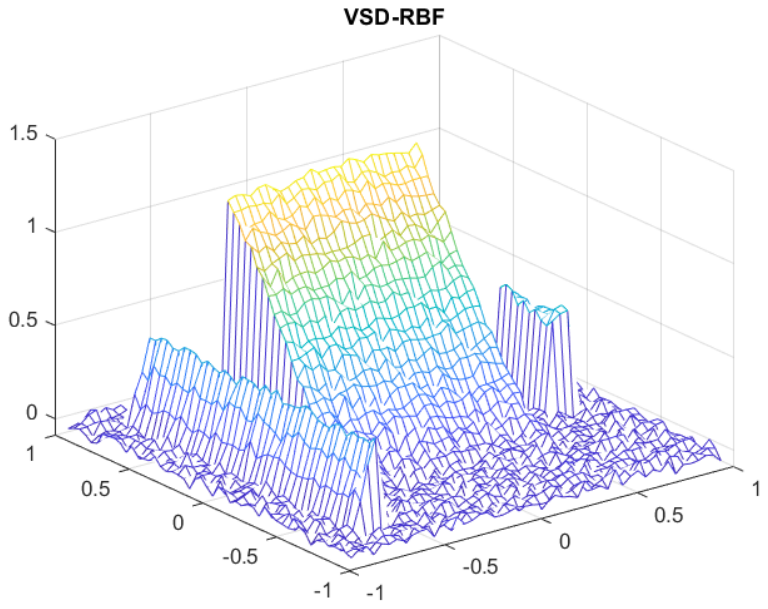}
        \includegraphics{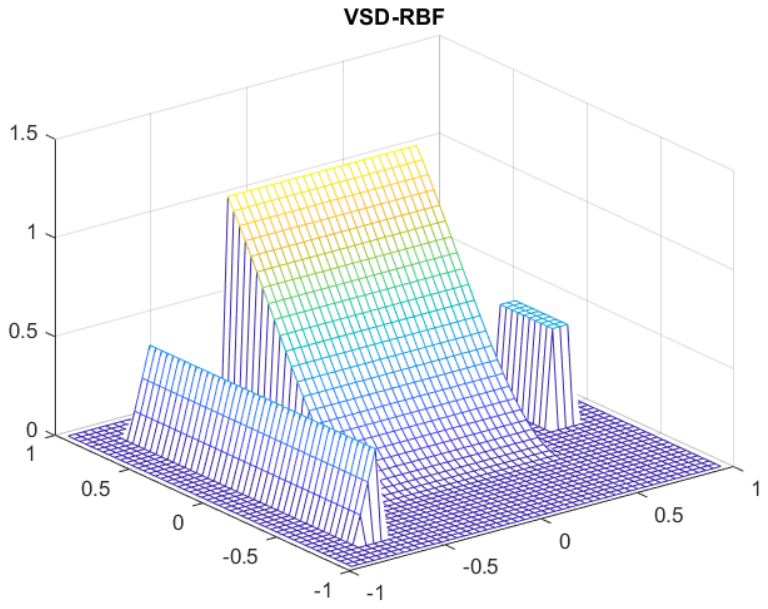}
        \includegraphics{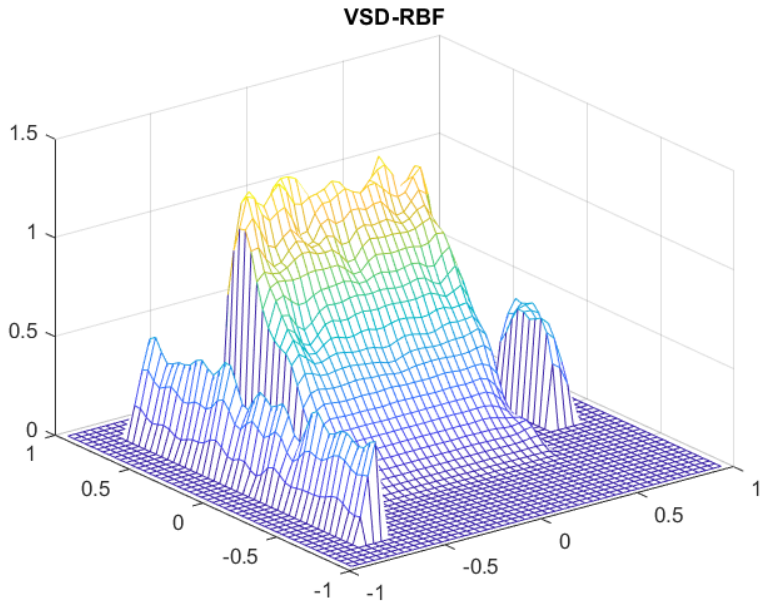}
        
        \includegraphics{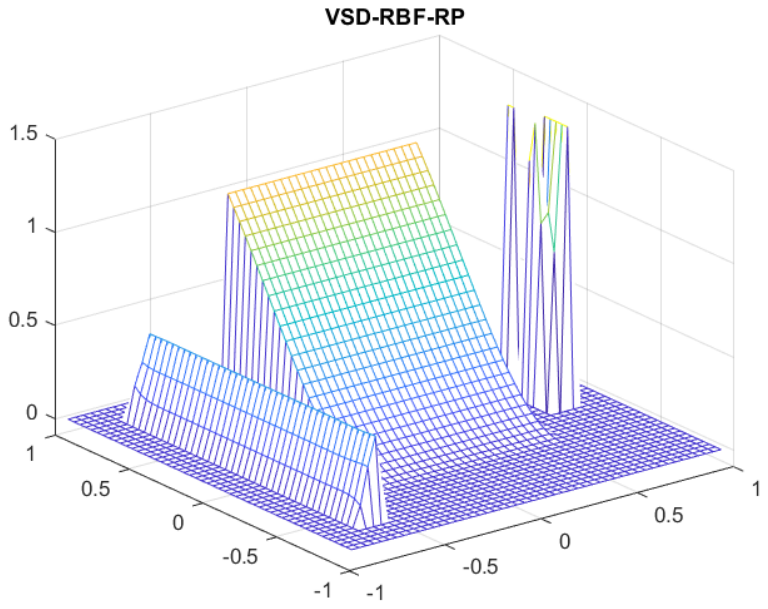}
        \includegraphics{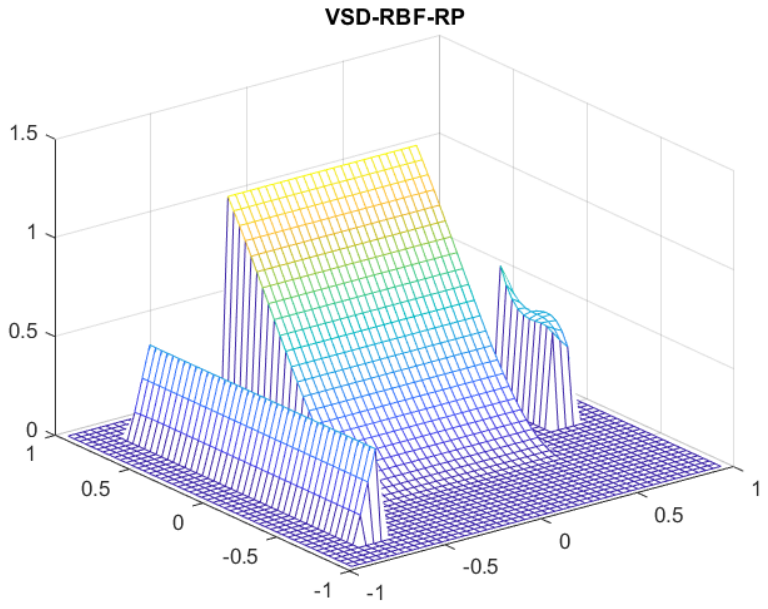}
        \includegraphics{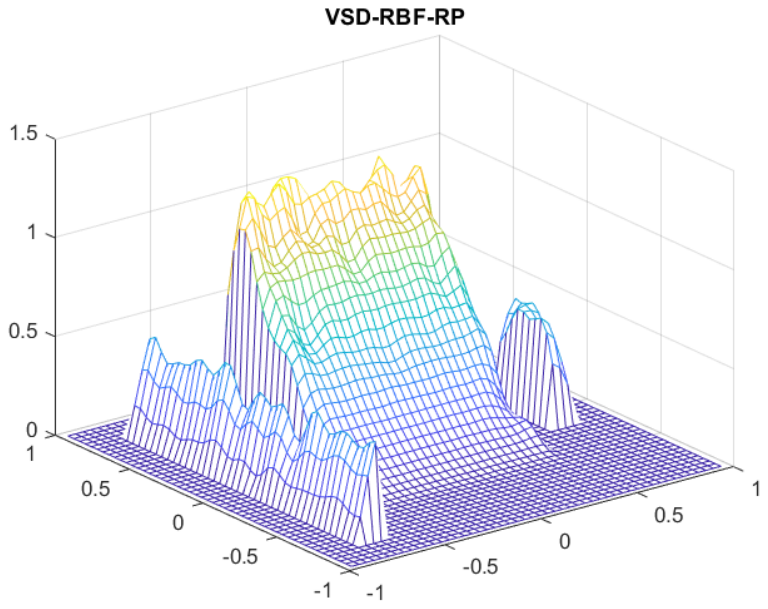}
	  \caption{\newm{Comparison of the approximate solutions for $f_2$ with $N = 625$ and $M = 2500$. Rows (top to bottom): RBF, RBF-RP, VSD-RBF, and VSD-RBF-RP \newn{methods}. Columns (left to right): $\varepsilon=0.1,1,$ and $10$.}}
\label{fig6}
\end{figure}
\begin{figure}
	\centering
        \setkeys{Gin}{width=7cm}
        \includegraphics{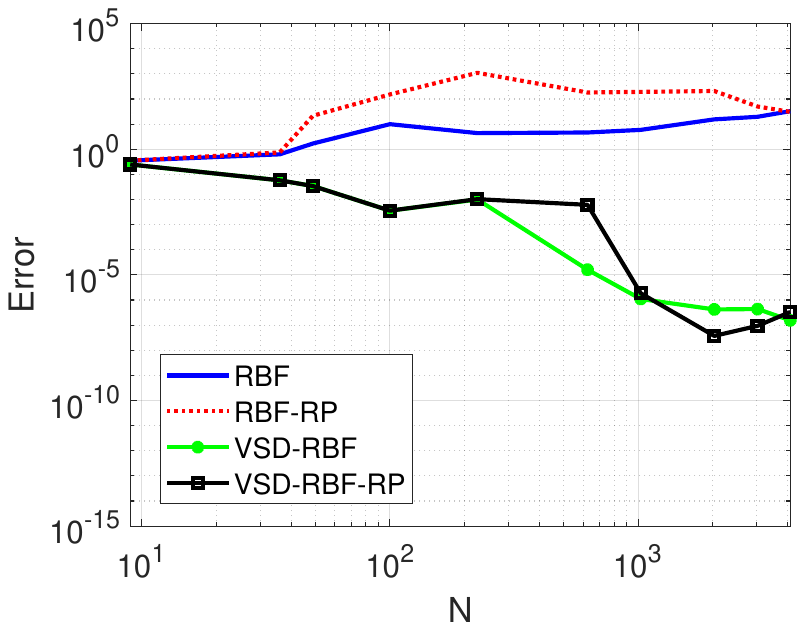}
	  \caption{\newm{Convergence behavior with respect to the number of centers $N$ at a fixed $\varepsilon = 1$.}}
\label{fig5_2}
\end{figure}
However, as shown in Figure~\ref{fig5_2}, the accuracy of VSD-RBF-RP \newn{interpolant} improves as the number of nodes $N$ increases, successfully overcoming the localized oscillations. To highlight this asymptotic advantage, the results for $N = 1600$ are shown in Figure~\ref{fig7}.
\begin{figure}
	\centering
        \setkeys{Gin}{width=7cm}
        \includegraphics{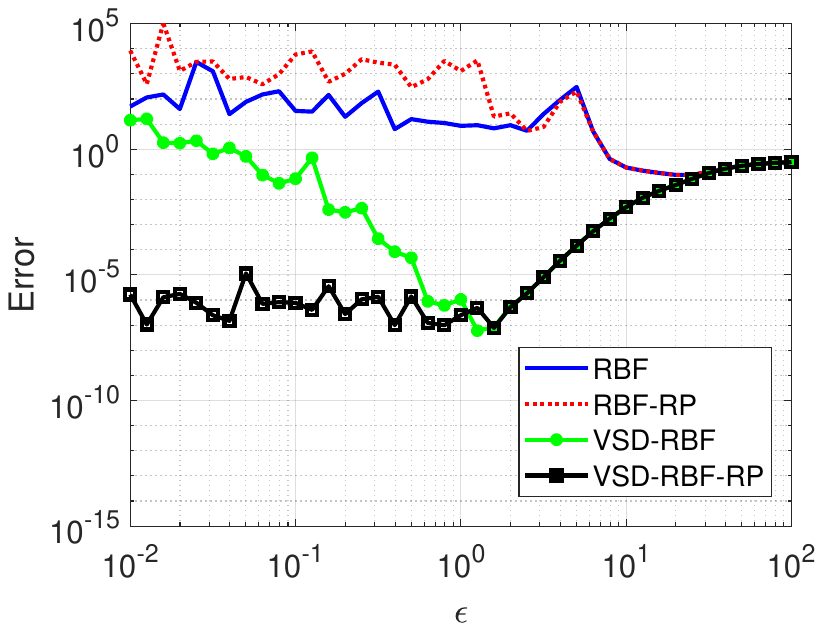}
	  \caption{\new{RMSE comparison for the function $f_2$ across various shape parameters and basis functions, with $N = 1600$ and $M = 6400$.}}
\label{fig7}
\end{figure}

Although VSD-RBF-RP \newn{method} clearly outperforms VSD-RBF \newn{method} at this higher number of centers, the structural similarity between $f_2$ and $f_1$ within their continuous subdomains raises \newm{the question of} why \newm{their RMSE} behaviors differ so substantially\newm{. To address this, Figure \ref{fig8}  displays the RMSE contour plots for the reconstruction of $f_2$ across different \newn{interpolation} bases for both $N=625$ and $N=1600$. These results demonstrate that the number of discontinuities does not inherently dictate the performance; instead, the discrepancy is predominantly governed by the geometry of the discontinuity boundaries, which exhibit sharp, rectangular corners in the case of $f_2$.}
\begin{figure}
	\centering
        \setkeys{Gin}{width=16cm,height=9cm,keepaspectratio}
        \includegraphics{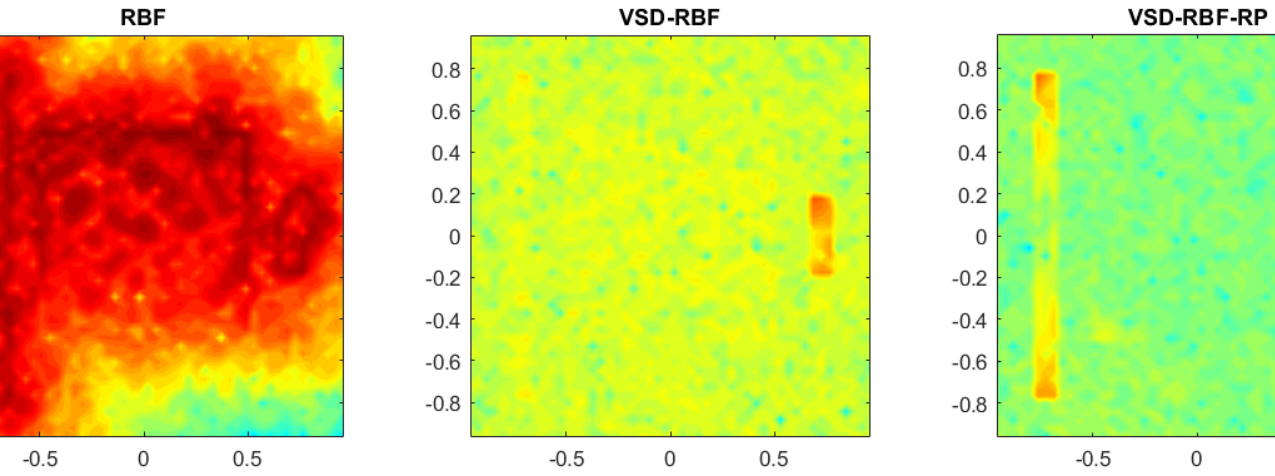}
        \includegraphics{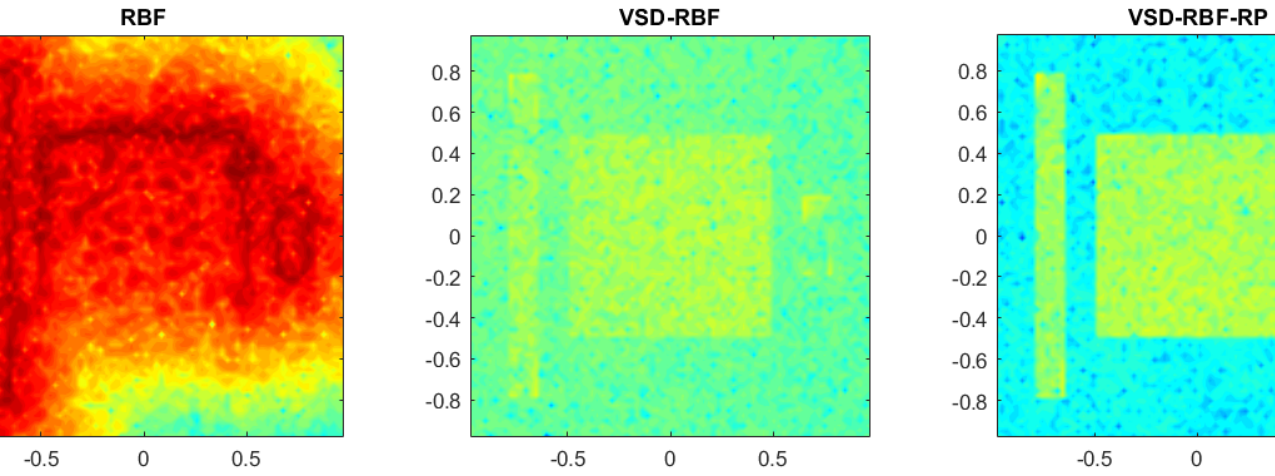}
	  \caption{\newm{Contour plots of the RMSE for the reconstruction of $f_2$ using different bases. For each method, the near-optimal shape parameters ($\varepsilon$) were determined through numerical experimentation. Top row ($N=625$, $M = 2500$): shows RBF ($\varepsilon=15.8489$), VSD-RBF ($\varepsilon=0.6310$), and VSD-RBF-RP ($\varepsilon=1.2589$) \newn{methods}. Bottom row ($N = 1600$, $M = 6400$):  RBF ($\varepsilon=19.9526$), VSD-RBF ($\varepsilon=1.2589$), and VSD-RBF-RP ($\varepsilon=1.5849$) \newn{methods}. In both experiments, the error remains notably elevated within the localized region $0.65 \leq x \leq 0.8$, $|y| \leq 0.2$.}}
\label{fig8}
\end{figure}
\end{example}
While Examples~\ref{exa1} and~\ref{exa2} demonstrate that the VSD-RBF-RP basis becomes increasingly effective compared to VSD-RBF \newn{basis} as the number of data points increases, numerical stability remains equally critical. To address this aspect, \newm{in} the following\newm{,} we compare the condition numbers obtained for the interpolation of $f_1$ and $f_2$. It is important to note that, while the condition number is not generally determined by the target function itself, in the present setting the interpolation nodes are generated through the specific definitions of $c_1(x)$ and $c_2(x)$. Consequently, the condition numbers associated with $f_1$ and $f_2$ may differ.
\begin{figure}
	\centering
        \setkeys{Gin}{width=7cm}
        \includegraphics{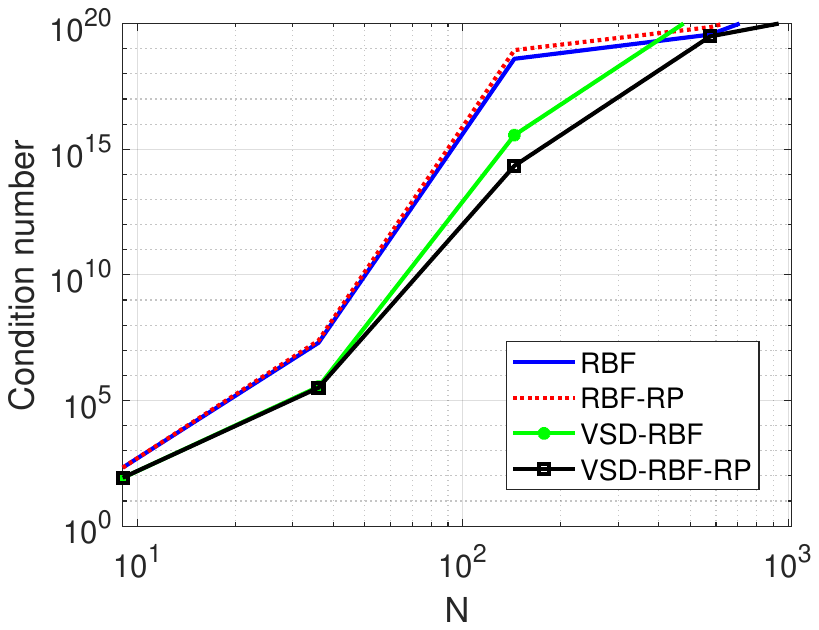}
        \includegraphics{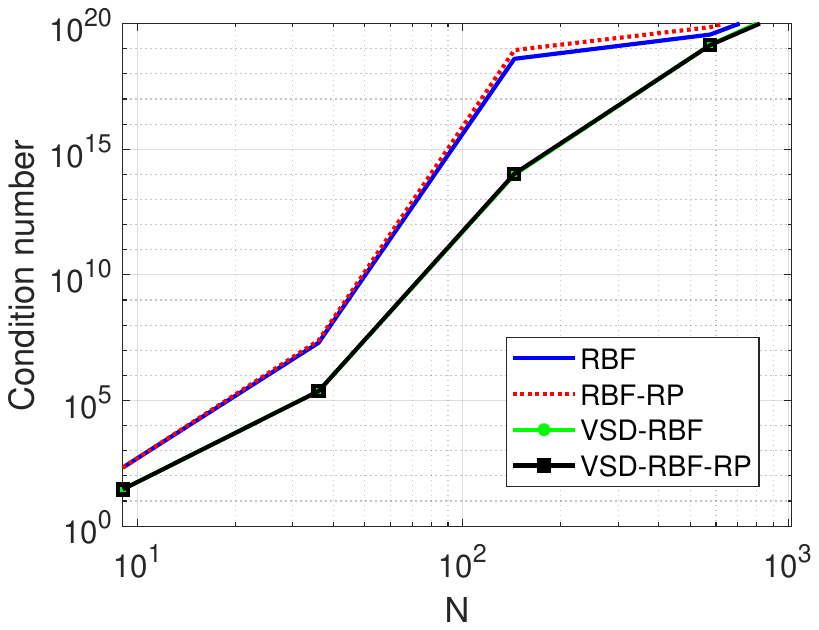}
	  \caption{Condition numbers of the coefficient matrices for the interpolation of $f_1$ and $f_2$ using the RBF \newm{($A$)}, RBF-RP \newm{($C$)}, VSD-RBF \newm{($A^c$)}, and VSD-RBF-RP \newm{($C^c$)} bases with shape parameter $\varepsilon = 1$. (from left to right).}
\label{fig9}
\end{figure} 

As shown in Figure \ref{fig9}, \newn{while the condition numbers for the standard RBF and RBF-RP methods remain identical for both $f_1$ and $f_2$, the VSD-based approaches exhibit distinct behaviors across the two test functions. Furthermore, as $N$ increases, the VSD-RBF-RP method demonstrates a moderately controlled growth in condition number compared to the other schemes, substantiating its numerical advantage in these interpolation configurations.}
\section{Applications in image processing} 
Interpolation plays an important role in several image-processing tasks, including image reconstruction, denoising, and resolution enhancement. In this section, we apply the proposed VSD-RBF-RP method to improve the clarity of phantom data. \newe{W}e evaluate the performance of the proposed methods using the bull's-eye and crescent-shaped phantoms introduced in \cite{DeMarchi16}. To define the corresponding scaling functions $c(x)$, we adopt the \newn{distance-based scaling strategy proposed in} \cite{DeMarchi20} \newn{(where $c(x)$ is set proportional to a unit distance of $1$).} A key objective of these tests is to demonstrate the capability of our methods in reconstructing high-resolution data from a limited number of sampling points. Specifically, we utilize $N=36$ clustered Halton points distributed in $[-1,1]^2$ as centers. To assess the approximation accuracy, a total of $M=14400$ evaluation points are employed in both test cases.
\begin{example}
We first consider the bull's-eye phantom, defined by the function $f_3$.
\begin{align*}
f_3(x,y)&=\begin{cases}
\frac{1}{2}, & x^2+y^2\leq \frac{1}{16}, \\
\frac{1}{4}, & \frac{1}{16}<x^2+y^2\leq \frac{1}{4}, \\
1, & \frac{1}{4}<x^2+y^2\leq \frac{9}{16},  \\
0, & \frac{9}{16}<x^2+y^2
\end{cases}
\end{align*}
\begin{align*}
c_3(x,y)&=\begin{cases}
1, & x^2+y^2\leq \frac{1}{16}, \\
2, & \frac{1}{16}<x^2+y^2\leq \frac{1}{4}, \\
3, & \frac{1}{4}<x^2+y^2\leq \frac{9}{16},  \\
4, & \frac{9}{16}<x^2+y^2
\end{cases}
\end{align*}
\begin{figure}
	\centering
        \setkeys{Gin}{width=7cm}
        \includegraphics{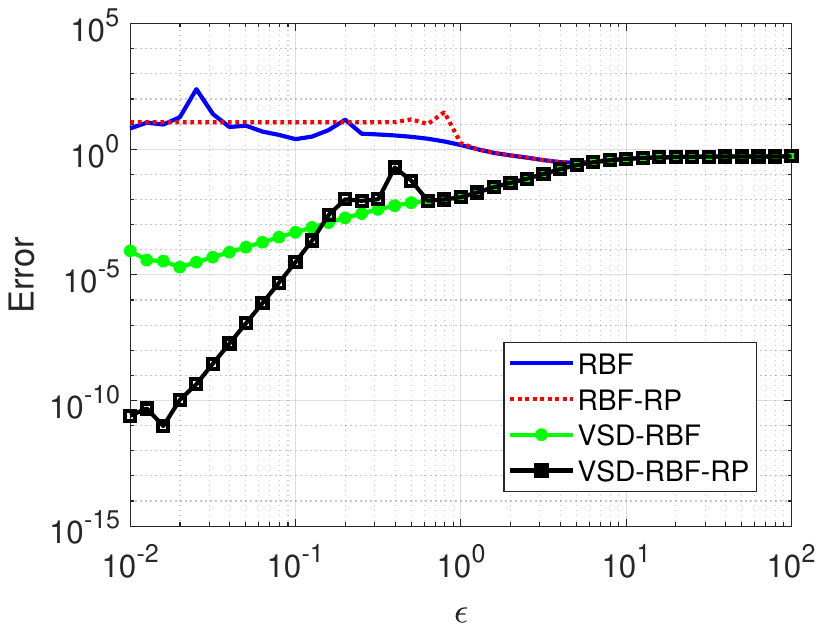}
        \includegraphics{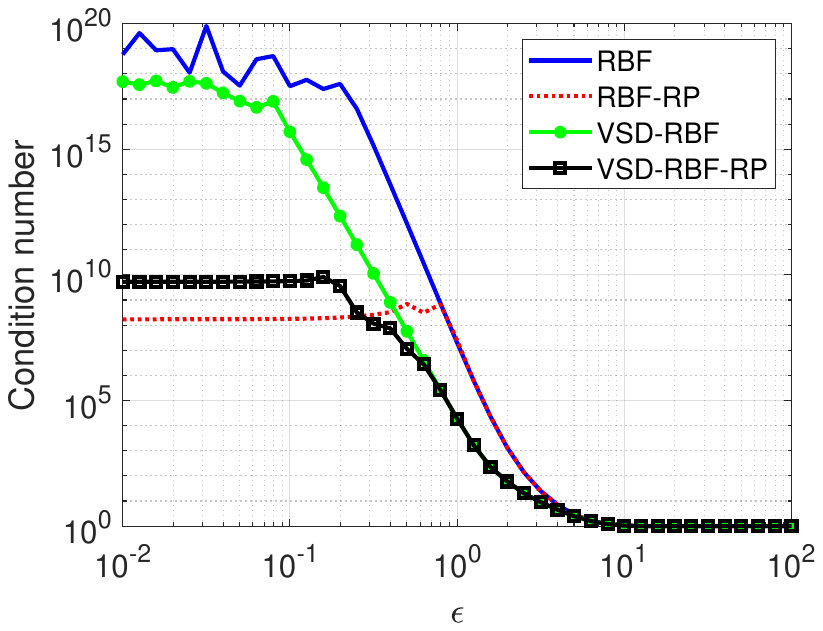}
	  \caption{\newm{RMSE and condition number for the test function $f_3$: (left) RMSE versus different shape parameters and basis functions for $N=36$ centers and $M=14400$ evaluation points; (right) the corresponding condition numbers.}}
\label{fig10}
\end{figure}
The accuracy and stability results are presented in Figure \ref{fig10}. As illustrated in this figure, the VSD-RBF-RP method not only achieves higher accuracy than VSD-RBF \newn{method} but also \new{exhibits more stable behavior.} Figure \ref{fig11} displays the reconstructed images obtained by applying different bases to interpolate $f_3$, generated using the \texttt{imshow()} command. 

\begin{figure}
	\centering
        \setkeys{Gin}{width=14cm,height=18cm,keepaspectratio}
        \includegraphics{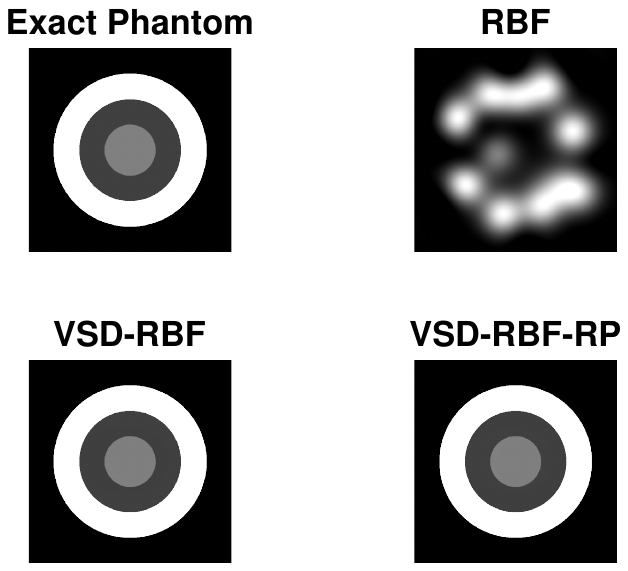}
	  \caption{\newm{Visual comparison of the reconstruction of $f_3$ on a $600 \times 600$ evaluation grid with a $6 \times 6$ center distribution ($N=36$): exact phantom, RBF ($\varepsilon=5.0119$), VSD-RBF ($\varepsilon=0.02$), and VSD-RBF-RP ($\varepsilon=0.0158$) \newn{methods}.}}
\label{fig11}
\end{figure}
It can be observed that the images produced by the VSD-RBF and VSD-RBF-RP bases closely resemble the original phantom data. Furthermore, the contour plots of the RMSE in Figure \ref{fig12} reveal that the VSD-RBF-RP method provides the highest level of accuracy among all tested bases.

\begin{figure}
	\centering
        \setkeys{Gin}{width=16cm,height=18cm,keepaspectratio}
        \includegraphics{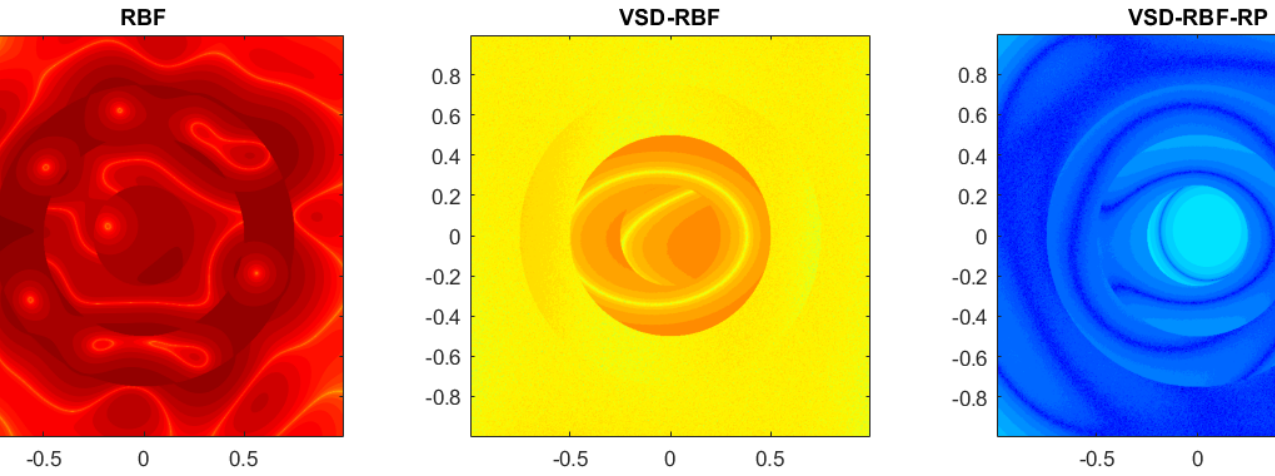}
	  \caption{\newm{Contour plots of RMSE for the reconstruction of $f_3$ evaluated on a $600 \times 600$ grid using $N=36$ centers with different bases. For each method, the near-optimal shape parameters ($\varepsilon$) were determined through numerical experimentation: RBF ($\varepsilon=5.0119$), VS-RBF ($\varepsilon=0.0200$), and VS-RBF-RP ($\varepsilon=0.0158$) \newn{methods}.}}
\label{fig12}
\end{figure}
\end{example}
\begin{example}
Next, we consider the crescent-shaped phantom defined by $f_4$. \newm{The accuracy and stability results, obtained under the same conditions as in Example 3,} are reported in Figure \ref{fig13}.
\begin{align*}
f_4(x,y)&=\begin{cases}
\frac{1}{2}, & (x-\frac{1}{8})^2+y^2\leq \frac{9}{64}, \\
1, & x^2+y^2\leq \frac{1}{4},(x-\frac{1}{8})^2+y^2> \frac{9}{64} , \\
0, & \frac{1}{4}<x^2+y^2
\end{cases}
\end{align*}
\begin{align*}
c_4(x,y)&=\begin{cases}
1, & (x-\frac{1}{8})^2+y^2\leq \frac{9}{64}, \\
2, & x^2+y^2\leq \frac{1}{4},(x-\frac{1}{8})^2+y^2> \frac{9}{64} , \\
3, & \frac{1}{4}<x^2+y^2
\end{cases}
\end{align*}
\begin{figure}
	\centering
        \setkeys{Gin}{width=7cm}
        \includegraphics{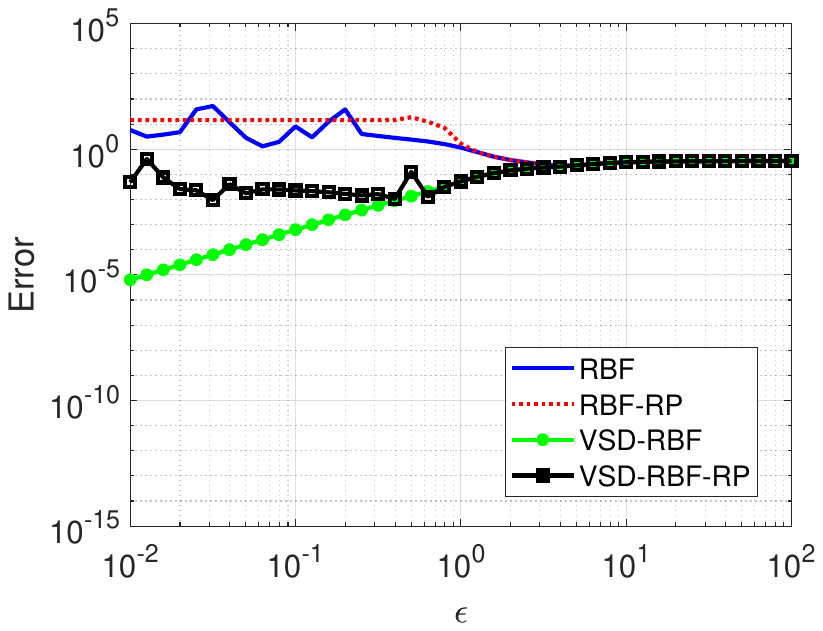}
        \includegraphics{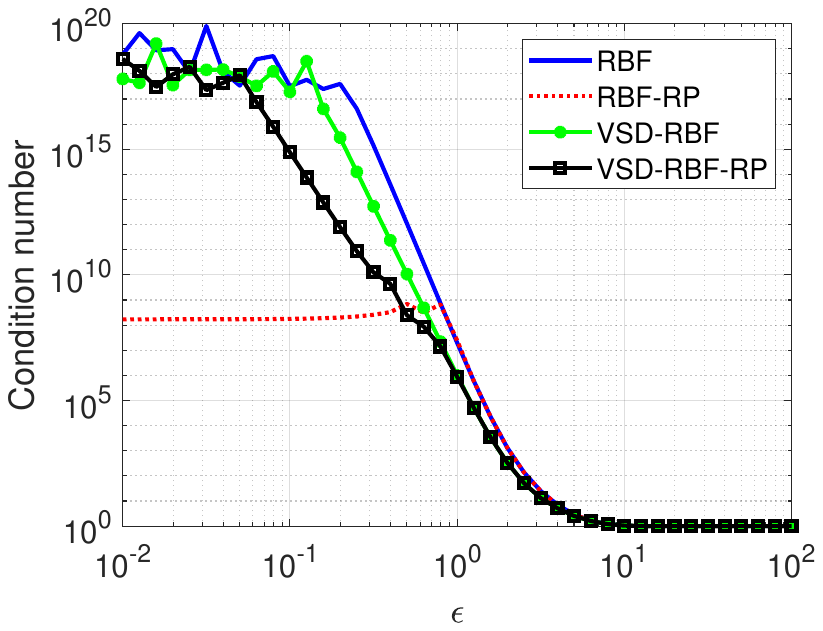}
	   \caption{\newm{RMSE and condition number for the test function $f_4$: (left) RMSE versus different shape parameters and basis functions for $N=36$ centers and $M=14400$ evaluation points; (right) the corresponding condition numbers.}}
\label{fig13}
\end{figure}
To mitigate the limited reconstruction accuracy observed in this setup, we increased the number of centers to $N=625$\newm{; the corresponding} stability and accuracy profiles are presented in Figure \ref{fig14}.
\begin{figure}
	\centering
        \setkeys{Gin}{width=7cm}
        \includegraphics{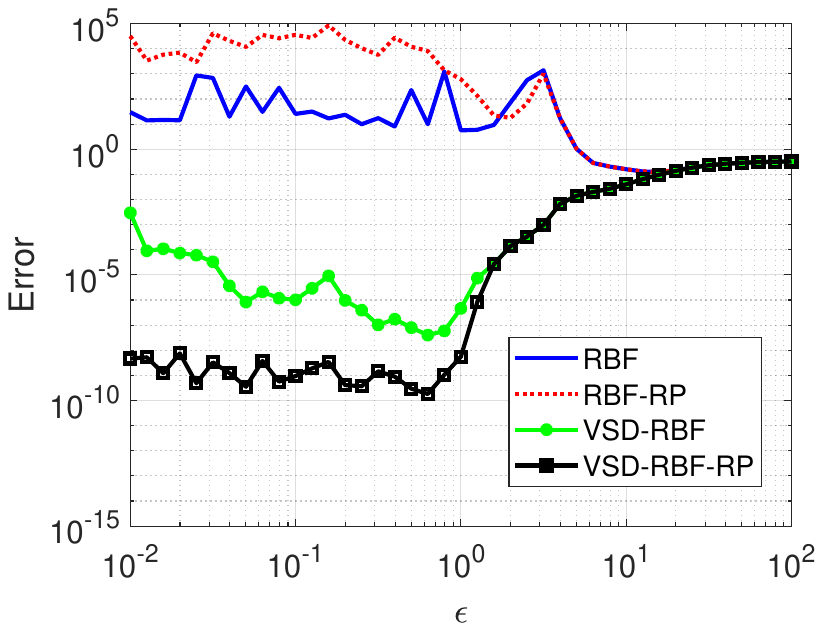}
        \includegraphics{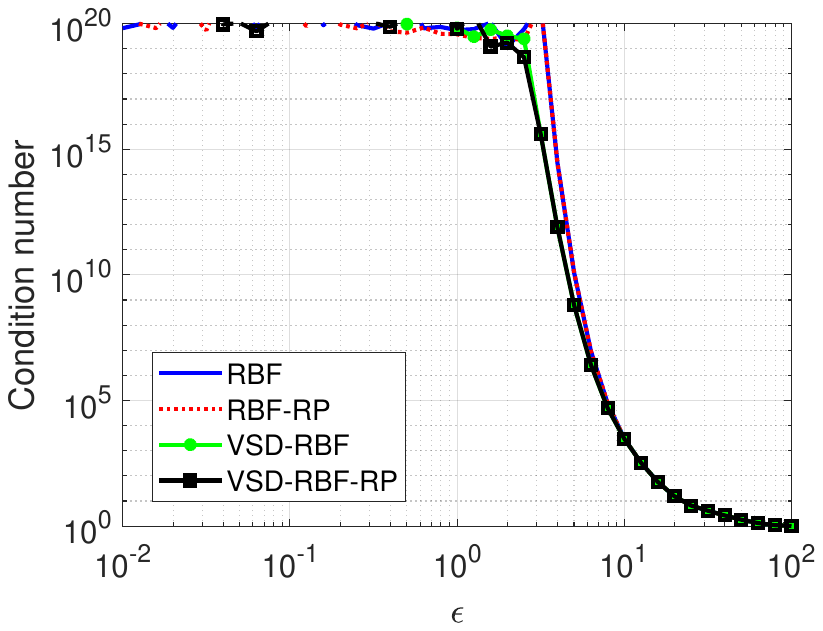}
	  \caption{\newm{RMSE and condition number for the test function $f_4$: (left) RMSE versus different shape parameters and basis functions for $N=625$ centers and $M=10000$ evaluation points; (right) the corresponding condition numbers.}}
\label{fig14}
\end{figure}
\new{As illustrated, the overall reconstruction performance improves significantly with this increased number of centers. However, it should be noted that \newn{the interpolation matrix becomes increasingly ill-conditioned as $N$ grows. \newe{Also,} non-singularity is not theoretically guaranteed for any $N$ \cite{Pooladi5}}.}

Consequently, \newm{the results for  $N=625$ are presented in Figures \ref{fig15} and \ref{fig16}. Figure \ref{fig15} displays the reconstructed images using different bases}, where both the VSD-RBF and VSD-RBF-RP methods produce visually clearer images.
\begin{figure}
	\centering
        \setkeys{Gin}{width=14cm,height=18cm,keepaspectratio}
        \includegraphics{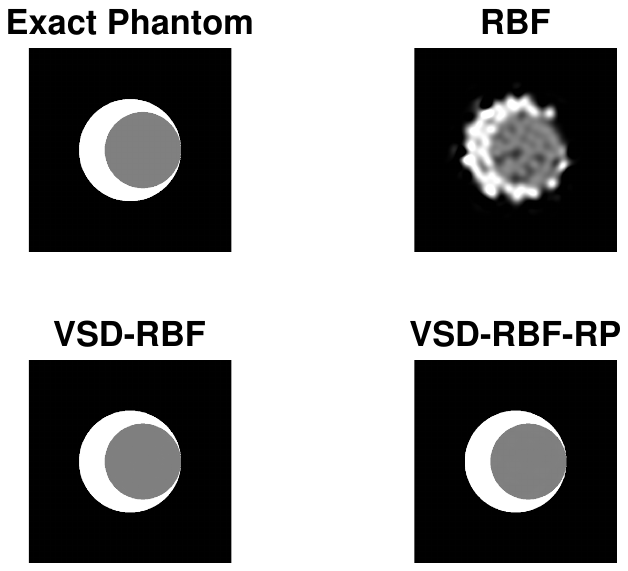}
	  \caption{\newm{Visual comparison of reconstruction results evaluated on a $600 \times 600$ evaluation grid with a $25 \times 25$ center distribution ($N=625$): Exact phantom, RBF ($\varepsilon=15.8489$), VSD-RBF ($\varepsilon=0.6310$), and VSD-RBF-RP ($\varepsilon=0.6310$) \newn{methods}.}}
\label{fig15}
\end{figure}

Furthermore, the RMSE contour plots in Figure \ref{fig16} confirm that the VSD-RBF-RP method yields the \newm{superior} reconstruction accuracy among all evaluated approaches, although higher error values \newm{remain observable} in the corner regions of the domain.
\begin{figure}
	\centering
        \setkeys{Gin}{width=16cm,height=18cm,keepaspectratio}
        \includegraphics{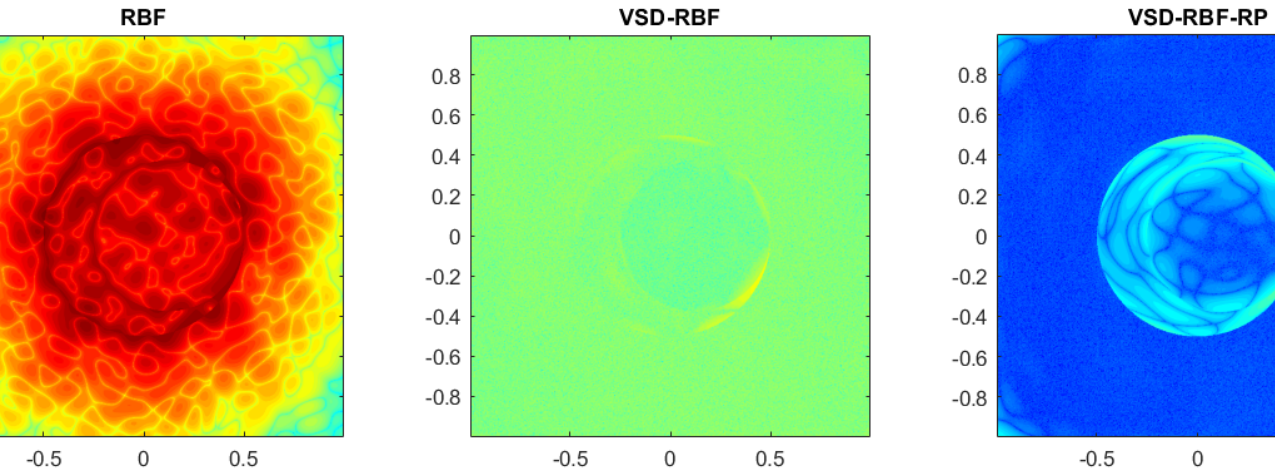}
	  \caption{\newm{Contour plots of RMSE for the reconstruction of $f_4$ evaluated on a $600 \times 600$ grid using $N=625$ centers with different bases. For each method, the near-optimal shape parameters ($\varepsilon$) were determined through numerical experimentation: RBF ($\varepsilon=15.8489$), VSD-RBF ($\varepsilon=0.6310$), and VSD-RBF-RP ($\varepsilon=0.6310$) \newn{methods}.}}
\label{fig16}
\end{figure}
\end{example}

\begin{figure}
	\centering
        \setkeys{Gin}{width=7cm}
        \includegraphics{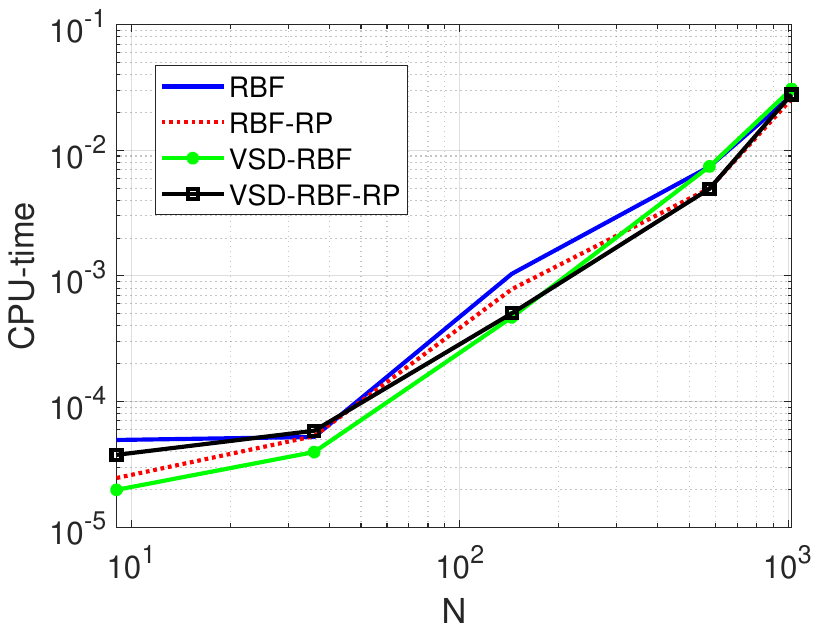}
        \includegraphics{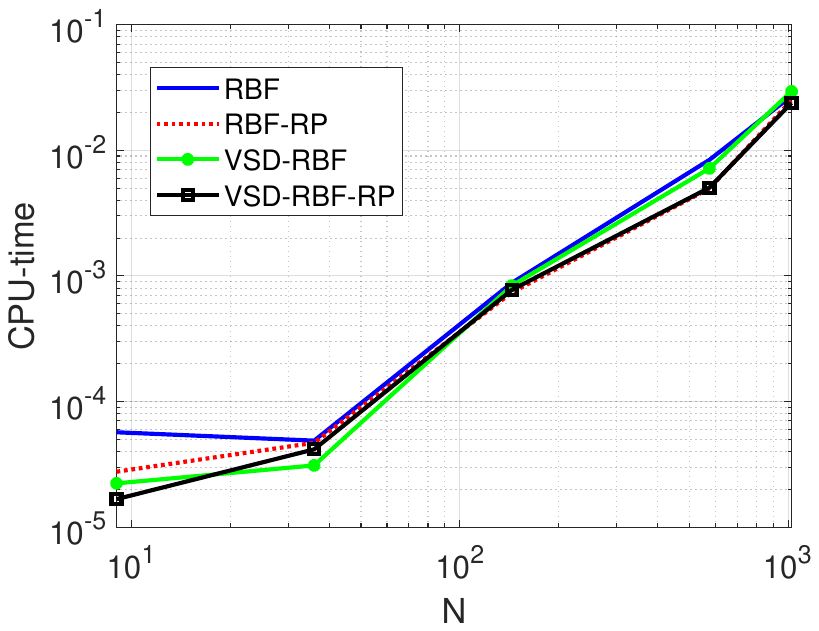}
	 \caption{\newm{Computational cost versus the number of centers ($N$) on a log--log scale for RBF ($A$), RBF-RP ($C$), VSD-RBF ($A^c$), and VSD-RBF-RP ($C^c$). The left and right panels correspond to $f_3$ and $f_4$, respectively, with $\varepsilon=1$ in both cases.}}
\label{fig17}
\end{figure}

The computational efficiency \newm{(measured in CPU-time)} of the investigated methods \newm{for reconstructing $f_3$ and $f_4$} is illustrated in Figure \ref{fig17}. The log-log plot demonstrates that all methods follow a similar asymptotic growth trend as the number of centers $N$ increases, which is consistent with the standard computational complexity of radial basis function interpolation. It is particularly noteworthy that the VSD-RBF-RP method exhibits a computational cost comparable to that of the classical RBF approach. This confirms that the enrichment of the basis functions with radial polynomials does not impose a significant additional computational burden, rendering the proposed VSD-RBF-RP method a computationally feasible choice for high-resolution approximations. It should be noted, however, that for smaller datasets where $N \leq 100$, the VSD-RBF method may offer a more cost-effective alternative.
\section{Conclusion and future work}
In this study, we introduced the VS-RBF-RP \newn{and} VSD-RBF-RP framework to enhance the performance and accuracy of classical VS-RBF \newn{and} VSD-RBF approximations. By incorporating degree-dependent radial polynomial terms into the lifted higher-dimensional formulation, the proposed approach significantly improves reconstruction \newm{accuracy}. Extensive numerical experiments validate that \newn{the} VSD-RBF-RP \newn{method achieves noticeably improved accuracy over VSD-RBF} across a broad spectrum of shape parameter choices.

Nevertheless, it is worth noting that the present formulation assumes the locations of the underlying discontinuities or internal interfaces are known a priori to construct the piecewise scale function $c(x)$. In practical applications where such geometric features are unknown beforehand, integrating our approach with automated discontinuity detection techniques ( such as gradient-based edge detectors, jump-detection algorithms, or machine learning-driven segmentation) represents a natural and promising direction for future research.
\section{Acknowledgments}
This work is based upon research funded by Iran National Science Foundation (INSF) under project number No.~403336.

\bibliographystyle{plain}
\bibliography{mybibfile}

@article{DeMarchi16,
  title={Kernel-based Image Reconstruction from Scattered {R}adon Data},
  author={S. De Marchi and A. Iske and A. Sironi},
  journal={Dolomites research notes on approximation},
  pages={1-12},
  year={2016},
  publisher={}
}

@article{DeMarchi20,
  title={Jumping with variably scaled discontinuous kernels ({VSDK}s)},
  author={S. De Marchi and F. Marchetti and E. Perracchione},
  journal={BIT Numerical Mathematics},
  pages={441-463},
  year={2020},
  publisher={}
}

@article{Bozzini14,
  title={Interpolation with variably scaled kernels},
  author={M. Bozzini and L. Lenarduzzi and M. Rossini and R. Schaback},
  journal={IMA Journal of Numerical Analysis},
  pages={1-21},
  year={2014},
  publisher={}
}

@article{Pooladi1,
  title={Radial polynomials as alternatives to smooth radial basis functions},
  author={F. Pooladi and H. Hosseinzadeh},
  journal={Journal of Mathematical Modeling},
  pages={337-354},
  year={2024},
  publisher={}
}

@article{Pooladi3,
  title={Stabilized interpolation using radial basis functions augmented with selected radial polynomials},
  author={F. Pooladi and E. Larsson},
  journal={Journal of Computational and Applied Mathematics},
  pages={115482},
  year={2024},
  publisher={}
}

@article{Romani19,
  title={Edge detection methods based on {RBF} interpolation},
  author={L. Romani and M. Rossini and D. Schenone},
  journal={Journal of Computational and Applied Mathematics},
  pages={532-547},
  year={2019},
  publisher={}
}

@article{Pooladi4,
  title={COMPARISON OF THREE SUITABLE APPROXIMATION SPACES FOR THE INTERPOLATION PROBLEM},
  author={F. Pooladi and E. Larsson},
  journal={The 10th Seminar on Numerical Analysis and its Applications},
  pages={},
  year={2024},
  publisher={}
}

@article{Pooladi5,
  title={Parameter Tuning in Hybrid {RBF-RP} Interpolation},
  author={F. Pooladi and M. Mohammadi and E. Larsson},
  journal={Engineering Analysis with Boundary Elements},
  year={2026},
  note={In press}
}

@article{Heidari23,
  title={Curvature-Based Characterization of Radial Basis Functions: Application to Interpolation},
  author={M. Heidari and M. Mohammadi and S. De Marchi},
  journal={Mathematical Modelling and Analysis},
  pages={415--433},
  year={2023},
  publisher={}
}

@article{Hardy90,
  title={Theory and applications of the multiquadric-biharmonic method. 20 years of discovery 1968–1988},
  author={R. L. Hardy},
  journal={Computers \& Mathematics with Applications},
  pages={163-208},
  year={1990},
  publisher={}
}

@article{Kansa901,
  title={Multiquadrics--A scattered data approximation scheme with applications to computational fluid-dynamics—I surface approximations and partial derivative estimates},
  author={E. J. Kansa},
  journal={Computers \& Mathematics with Applications},
  pages={127-145},
  year={1990},
  publisher={}
}

@article{Kansa902,
  title={Multiquadric--A scattered data approximation scheme with applications to computational fluid-dynamics—II solutions to parabolic, hyperbolic and elliptic partial differential equations},
  author={E. J. Kansa},
  journal={Computers \& Mathematics with Applications},
  pages={147-161},
  year={1990},
  publisher={}
}

@article{Kansa923,
  title={A strictly conservative spatial approximation scheme for the governing engineering and physics equations over irregular regions and inhomogeneously scattered nodes},
  author={E. J. Kansa},
  journal={Computers \& Mathematics With Applications},
  pages={147-161},
  year={1992},
  publisher={}
}

@article{Kansa92,
  title={Improved accuracy of multiquadric interpolation using variable shape parameters},
  author={E.J. Kansa and R.E. Carlson},
  journal={Computers \& Mathematics with Applications},
  pages={99-120},
  year={1992},
  publisher={}
}

@article{Bozzini02,
  title={Adaptive interpolation by scaled multiquadrics},
  author={M. Bozzini and L. Lenarduzzi and R. Schaback},
  journal={Advances Computational Mathematics},
  pages={375-387},
  year={2002},
  publisher={}
}

@article{Driscoll07,
  title={Adaptive residual subsampling methods for radial basis function interpolation and collocation problems},
  author={T.A. Driscoll and A. Heryudono},
  journal={Computers \& Mathematics with Applications},
  pages={927-939},
  year={2007},
  publisher={}
}

@article{Fornberg07,
  title={The Runge phenomenon and spatially variable shape parameters in {RBF} interpolation},
  author={B. Fornberg and J. Zuev},
  journal={Computers \& Mathematics with Applications},
  pages={379-398},
  year={2007},
  publisher={}
}

@book{Wendland04,
  title={Scattered data approximation},
  author={Wendland, Holger},
  volume={17},
  year={2005},
  publisher={Cambridge university press, Cambridge}
}

@article{Fasshauer07,
  title={MESHFREE APPROXIMATION METHODS WITH MATLAB},
  author={G. Fasshauer},
  journal={},
  volume={},
  year={2007},
  publisher={illinois institute of technology}
}

@article{Cavoretto24,
  title={Comparing Deterministic and Statistical Optimization Techniques for the Shape Parameter Selection in {RBF} Interpolation},
  author={R. Cavoretto and A. De Rossi and A. Haider and S. Lancellotti},
  journal={Dolomites research notes on approximation},
  pages={48--55},
  year={2024},
}

@article{Sun23,
  title={An Adaptive Selection Method for Shape Parameters in {MQ-RBF} Interpolation for Two-Dimensional Scattered Data and Its Application to Integral Equation Solving},
  author={J. Sun and L. Wang and D. Gong},
  journal={ Fractal and Fractional},
  pages={488},
  year={2023},
}

@article{Audone25,
  title={A recipe for learning Variably Scaled Kernels via Discontinuous Neural Networks},
  author={G. Audone and F. Della Santa and E. Perracchione and S. Pieraccini},
  journal={Journal of Computational and Applied Mathematics},
  pages={116653},
  year={2025},
}

@article{Larsson24,
  title={Scaling of radial basis functions},
  author={E. Larsson and R. Schaback},
  journal={IMA Journal of Numerical Analysis},
  pages={1130-1152},
  year={2024},
}

@article{Mohammadi26,
  title={Spatially variable shape parameters for {RBF} interpolation: A curvature-based adaptive strategy},
  author={M. Mohammadi and M. Heidari},
  journal={Engineering Analysis with Boundary Elements},
  pages={106923},
  year={2026},
}

@article{Cavoretto21,
  title={On the search of the shape parameter in radial basis functions using univariate global optimization methods},
  author={R. Cavoretto and A. De Rossi and M. S Mukhametzhanov and Y. D Sergeyev},
  journal={Journal of Global Optimization},
  pages={305-327},
  year={2021},
}

@article{Esmaeilbeigi14,
  title={A new approach based on the genetic algorithm for finding a good shape parameter in solving partial differential equations by {K}ansa's method},
  author={M. Esmaeilbeigi and MM. Hosseini},
  journal={Applied mathematics and computation},
  pages={419-428},
  year={2014},
}

@article{Chen23,
  title={On the selection of a better radial basis function and its shape parameter in interpolation problems},
  author={C. S. Chen and A. Noorizadegan and D. L. Young and CS Chen},
  journal={Applied Mathematics and Computation},
  pages={127713},
  year={2023},
}

@article{Fornberg08,
  title={On choosing a radial basis function and a shape parameter when solving a convective {PDE} on a sphere},
  author={B. Fornberg and C. Piret},
  journal={Journal of Computational Physics},
  pages={2758-2780},
  year={2008},
}



\end{document}